\documentclass[11pt]{amsart}
\usepackage{amsxtra,amscd, color}
\usepackage{amssymb}
\theoremstyle{plain}

\newcommand{\cleqn}{\setcounter{equation}{0}}
\newcommand{\clth}{\setcounter{theorem}{0}}
\newcommand {\sectionnew}[1]{\section{#1}\cleqn\clth}
\newcommand{\nn}{\hfill\nonumber}
\newtheorem{theorem}{Theorem}[section]
\newtheorem{lemma}[theorem]{Lemma}
\newtheorem{definition-theorem}[theorem]{Definition-Theorem}
\newtheorem{proposition}[theorem]{Proposition}

\newtheorem{question}[theorem]{Question}
\newtheorem{corollary}[theorem]{Corollary}
\newtheorem{definition}[theorem]{Definition}
\newtheorem{example}[theorem]{Example}

\newtheorem{remark}[theorem]{Remark}
\newtheorem{remarks}[theorem]{Remarks}
\newtheorem{conjecture}[theorem]{Conjecture}

\newcommand \bth[1] { \begin{theorem}\label{t#1} }
\newcommand \ble[1] { \begin{lemma}\label{l#1} }

\newcommand \bpr[1] { \begin{proposition}\label{p#1} }
\newcommand \bqu[1] { \begin{question}\label{q#1} }
\newcommand \bco[1] { \begin{corollary}\label{c#1} }
\newcommand \bde[1] { \begin{definition}\label{d#1}\rm }
\newcommand \bex[1] { \begin{example}\label{e#1}\rm }
\newcommand \bre[1] { \begin{remark}\label{r#1}\rm }
\newcommand \bres[1] { \begin{remarks}\label{r#1}\rm }
\newcommand \bcj[1] { \begin{conjecture}\label{j#1}\rm }
\renewcommand {\eth} { \end{theorem} }
\newcommand {\ele} { \end{lemma} }

\newcommand {\epr} { \end{proposition} }
\newcommand {\equ} {\end{question} }
\newcommand {\eco} { \end{corollary} }
\newcommand {\ede} { \end{definition} }
\newcommand {\eex} { \end{example} }
\newcommand {\ere} { \end{remark} }
\newcommand {\eres} { \end{remarks} }
\newcommand {\ecj} { \end{conjecture} }
\newcommand {\enota} { \end{notation} }
\newcommand \thref[1]{Theorem \ref{t#1}}
\newcommand \leref[1]{Lemma \ref{l#1}}
\newcommand \prref[1]{Proposition \ref{p#1}}

\newcommand \coref[1]{Corollary \ref{c#1}}

\newcommand \deref[1]{Definition \ref{d#1}}

\newcommand \reref[1]{Remark \ref{r#1}}

\def \Cset {{\mathbb C}}
\def \Aset {{\mathbb A}}
\def \KK {{\mathbb K}}
\def \Zset {{\mathbb Z}}
\def \Nset {{\mathbb N}}

\def \AA  {{\mathcal{A}}}           

\def \FF {{\mathcal{F}}}

\def \VV {{\mathcal{V}}}

\def \ZZ {{\mathcal{Z}}}

\def \ZZ {{\mathcal{Z}}}

\def \de {\delta}

\def \la {\lambda}

\def \de {\delta}

\def \sig {\sigma}

\def \ep {\epsilon}
\def \sig{\sigma}

\def \mt  {\mapsto}
\def \lra {\longrightarrow}

\def \hra {\hookrightarrow}

\def \rcor {\rangle}
\def \lcor {\langle}

\def \ol {\overline}

\def \wh {\widehat}

\def \id { {\mathrm{id}} }

\def \n  {\mathfrak{n}}
\def \mm  {\mathfrak{m}}

\DeclareMathOperator \Der { {\mathrm{Der}} }

\DeclareMathOperator \charr { {\mathrm{char}} }

\DeclareMathOperator \maxspec { {\mathrm{Maxspec}}}
\DeclareMathOperator \ord { {\mathrm{ord}} }

\DeclareMathOperator \Ker { {\mathrm{Ker}} }
\DeclareMathOperator \tr { {\mathrm{tr}} }

\DeclareMathOperator \End { {\mathrm{End}} }

\DeclareMathOperator \st   { {\mathrm{st}} }

\DeclareMathOperator \Irr { {\mathrm{Irr}} }

\DeclareMathOperator \red  { {\mathrm{red}}}
\DeclareMathOperator \reg  { {\mathrm{reg}}}

\renewcommand \Im { {\mathrm{Im}} }

\begin{document}
\title[Azumaya loci and discriminant ideals]
{Azumaya loci and discriminant ideals of PI algebras}
\author[K. A. Brown]{K. A. Brown}
\address{School of Mathematics and Statistics \\
University of Glasgow \\
Glasgow G12 8QW, Scotland
}
\email{Ken.Brown@glasgow.ac.uk}
\author[M. T. Yakimov]{M. T. Yakimov}
\thanks{The research of MY was supported by
 NSF grant DMS-1601862.}
\address{
Department of Mathematics \\
Louisiana State University \\
Baton Rouge, LA 70803 \\
U.S.A.
}
\email{yakimov@math.lsu.edu}
\date{}
\begin{abstract} We prove that, under mild assumptions,  for all positive integers $\ell$, the zero set of the discriminant ideal $D_{\ell}(R/Z(R), \tr)$ of a prime affine polynomial identity (PI) algebra $R$ coincides with the zero set of the modified discriminant ideal  $MD_{\ell}(R/Z(R), \tr)$ of $R$, and give an explicit description
of this set in terms of the dimensions of the irreducible representations of $R$.
Furthermore, we prove that, when $\ell$ is the square of the PI-degree of $R$,  this zero set is precisely the complement
of the Azumaya locus of $R$. This description is used to classify the Azumaya loci of the mutiparameter quantized Weyl algebras at roots
of unity. As another application, we prove that the zero set of the top discriminant ideal of a prime affine PI algebra $R$ coincides
with the singular locus of the center of $R$, provided that the discriminant ideal has height at least 2, $R$ has finite global dimension and $R$ is a Cohen--Macaulay module over its center.
\end{abstract}
\maketitle
\sectionnew{Introduction}
\label{intro}
\subsection{Setting}
\label{1.1}
Let $\KK$ be an algebraically closed field (of arbitrary characteristic), 
and let $R$ be a prime $\KK$-algebra which is a finitely generated module over its affine center $Z(R)$.
A fundamental invariant of $R$ is its {\em{Azumaya locus}}
\[
\AA(R) \subseteq \maxspec Z(R),
\]
the dense open subset of $\maxspec Z(R)$ which parametrizes the irreducible $R$-modules of maximal dimension. More precisely, $\mm \in \AA(R)$ if and only if
$\mm R$ is the annihilator in $R$ of an irreducible $R$-module $V$ with $\dim_{\KK}V = n$, where $n$ is the PI-degree of $R$, the maximal $\KK$-dimension of
irreducible $R$-modules.

The purpose of this paper is to show that, under mild hypotheses on the algebra $R$, the Azumaya locus is the complement of
the zero set of the top \emph{discriminant ideal}, that of order equal to the square of the PI-degree of $R$.
We also prove similar descriptions of the zero loci of the lower discriminant ideals of $R$ in terms of the irreducible representations 
of $R$ with a given central annihilator.

\emph{Discriminant ideals} are classical invariants which have been the focus of some recent work on automorphisms of PI-algebras, \cite{CPWZ, NTY}.
They are ideals of $Z(R)$ defined using
a \emph{trace map} from $R$ to $Z(R)$. The relevant definitions and examples are recalled in $\S\S$\ref{2.1}, \ref{2.1extra} for trace maps, and in $\S\S$\ref{2.1A} and \ref{discex} for discriminants.
We impose mild conditions on the definition of a trace map $\tr: R \to Z(R)$, namely that it is a cyclic, $Z(R)$-linear, non-zero map. 
In particular, this permits $\tr(1)$ to be $0$. For instance, the reduced trace on $R$ is a trace map in this sense for fields $\KK$ of arbitrary 
characteristic, as we shall show in \prref{re}.
\subsection{Results}
\label{1.2}
For an ideal $I$ of $Z(R)$, we will denote by $\VV(I)$ its zero set in $\maxspec Z(R)$. For $\mm \in \maxspec Z(R)$,
denote by $\Irr_\mm(R)$ the equivalence classes of finite dimensional irreducible representations of $R$
which are annihilated by $\mm$. Since $R$ is a finitely generated module over $Z(R)$, $R/\mm R$ is finite dimensional, and thus, 
$\Irr_\mm(R)$ is a finite set. By $Z(R)$-linearity, each trace $\tr : R \to Z(R)$ descends to a $\KK$-linear map 
\begin{equation}
\label{trnote}
\tr_\mm : R/ \mm R \to Z(R)/ \mm
\end{equation}
which is a trace on $R/\mm R$ provided it is nonzero. A trace
$\tr : R \to Z(R)$ will be called {\em{almost representation theoretic}} if, for each $\mm \in \maxspec Z(R)$, there exists a function
\begin{equation}
\label{rep-theor-tr1}
s_\mm : \Irr_\mm(R) \to \KK
\end{equation}
such that
\begin{equation}
\label{rep-theor-tr2}
\tr_\mm = \sum_{ V \in \Irr_\mm(R)} s_\mm(V) \tr_V
\end{equation}
where the traces in the right hand side are the standard trace maps on endomorphisms of finite dimensional vector spaces.
The function $s_\mm$ is allowed to take on zero values for some representations $V \in \Irr_\mm(R)$. When $s_{\mm}(V)$ is required to be non-zero at least for those $\mm$ in the Azumaya locus, we shall say that $\tr$ is \emph{representation theoretic}.
\medskip
\\
\noindent
{\bf{Main Theorem.}} {\em{Let $R$ be a prime affine algebra over an algebraically closed field $\KK$, which is
a finitely generated module over its center $Z(R)$. Assume that $R$ coincides with its trace ring $T(R)$; in particular this holds if $Z(R)$ is normal.
Let $n$ be the PI degree of $R$.

{\em{(}}a{\em{)}} Let $\tr$ be a representation theoretic trace map on $R$. The zero sets of the $n^2$-discriminant 
and modified $n^2$-discriminant ideals coincide and equal
the complement of the Azumaya locus of $R$. Namely,
\[
\VV(D_{n^2}(R/Z(R),\tr)) = \VV(MD_{n^2}(R/Z(R),\tr)) = \maxspec( Z(R))\; \backslash \;\AA(R).
\]

{\em{(}}b{\em{)}} The reduced trace on $R$ is representation theoretic.

{\em{(}}c{\em{)}} The standard trace on $R$ is representation theoretic provided that  $\charr \KK \notin [1,n]$.

{\em{(}}d{\em{)}} For all positive integers $\ell$ and any almost representation theoretic trace $\tr$,
the zero sets of the $\ell$-discriminant ideal and the modified $\ell$-discriminant ideal of $R$ coincide and are given in
terms of the irreducible representations of $R$ by
\begin{multline}\label{eq-master}
 \VV(D_{\ell}(R/Z(R),\tr) )  = \VV(MD_{\ell}(R/Z(R),\tr))
\\
=  \Big\{ \mm \in \maxspec Z(R) \mid
\sum_{V \in \Irr_\mm(R), s_\mm(V) \neq 0} (\dim_\KK V)^2  < \ell \Big\}.
\end{multline}

{\em{(}}e{\em{)}} Assume that  $\charr \KK \notin [1,n]$.  For all positive integers $\ell$ and for the reduced trace $\tr_{\red}$,
\begin{align}
\label{eq-partb}
\VV(D_{\ell}(R/Z(R),\tr_{\red}) )  &= \VV(MD_{\ell}(R/Z(R),\tr_{\red}))
\\
&=  \Big\{ \mm \in \maxspec Z(R) \mid
\dim_\KK \big( (R/\mm R)/J(R/\mm R)  \big) < \ell \Big\}
\nn
\\
&=\Big\{ \mm \in \maxspec Z(R) \mid
\sum_{V \in \Irr_\mm(R)} (\dim_\KK V)^2  < \ell \Big\}.
\nn
\end{align}. 

{\em{(}}f{\em{)}} For all traces on $R$ and for all integers $\ell$ with $\ell >n^2$,
\[
D_{\ell}(R/Z(R),\tr) = MD_{\ell}(R/Z(R),\tr) =0.
\]
}}
\smallskip

In \thref{red-tr-main} we prove a stronger version of part (e) of the theorem without any assumptions on the 
characteristic of the field $\KK$. 
The above result is particularly suited for reduction mod $p$ methods because the characteristic of the base
field can be arbitrary. 
In \thref{tCH} we prove a full analog of this theorem for all Cayley-Hamilton algebras
in the sense of Procesi \cite{P}, but that version requires $\charr \KK =0$ due to its dependence on the results in \cite{P}.

Denote by $J(A)$ the Jacobson radical of an algebra $A$. In \prref{MTa} we obtain a more general statement than that of 
part (a) of the Main Theorem for all traces $\tr \colon R \to Z(R)$ that have the property that they
descend to traces on the semisimple quotients $(R/\mm R)/J(R/\mm R)$ for all $\mm \in \maxspec Z(R)$.
(In particular, this fact is applicable to all almost representation theoretic traces.)
This descent condition explicitly means that
\begin{equation}\label{induce}
a \in R, \mm \in \maxspec Z(R),\;  a + \mm R \in J(R/\mm R) \; \;  \Rightarrow \; \; \tr (a) \in \mm.
\end{equation}

The theorem gives a new method to determine the Azumaya loci of PI algebras by computing their discriminant ideals.
We use this idea to determine the Azumaya loci of the multiparameter quantized Weyl algebras.
In the opposite direction, the theorem provides valuable information about the discriminant and modified discriminant ideals of a PI algebra whose
Azumaya locus is known. For one such application we refer to \cite{WWY}
where Poisson geometric and representation theoretic techniques
are used to determine the Azumaya loci of the 3 dimensional PI Sklyanin algebras and then our theorem can be applied to obtain information on the
discriminant ideals of these algebras. Those discriminant ideals are impossible to compute directly due to the lack of PBW bases of
the algebras in question.

In a third direction we relate the zero locus of the discriminant ideal of a prime PI algebra $R$ with the
singlular locus of its center $Z(R)$. Several results are proved, the strongest of which is that the discriminant ideal of $R$ coincides
with the singular locus of $Z(R)$, assuming that the former has height at least 2, $R$ has finite global dimension and
is a Cohen--Macaulay module over its center.

The Main Theorem is proved in Sect. \ref{pf}. Sect. \ref{disc} contains background on traces, discriminants and discriminant ideals.
Sect. \ref{CH} contains the statement and proof of an analog of the Main Theorem for Cayley-Hamilton algebras.
Sect. \ref{consequences} contains applications of the Main Theorem to the above mentioned relation between the discriminant ideal of a PI algebra $R$ and
the singular locus of $Z(R)$. Sect. \ref{qWeyl} addresses the second application of the Main Theorem that determines the Azumaya loci of the cocycle
twists of arbitrary tensor products of quantized Weyl algebras.

In the important case when $R$ is a Poisson order \cite{BrGo},
the Main Theorem establishes a bridge between the previous results on the Poisson properties of discriminants and Azumaya loci \cite{BrGo,NTY}.
This is described in Sect. \ref{pf}. This section also contains results that establish the needed normality of $Z(R)$ in the Main Theorem in
various common situations.
\subsection{Proof of part (a) of the Main Theorem.}
\label{MTproof}
Here is a sketch of the proof of part (a) of the Main Theorem.
Keep its notation and hypotheses, and let $\mm$ be a maximal ideal of $Z(R)$. Since $\tr$
is $Z(R)$-linear, it induces a trace-like map $\tr_{\mm}$, (which might be $0$), on $R/\mm R$, taking values in $\KK$, and hence also
a symmetric bilinear form  $\lcor -, - \rcor_{\tr}$ on $R/\mm$R, given by $\lcor a, b \rcor_{\tr} = \tr_{ \mm}(ab)$. Now
suppose that $\mm$ is in the Azumaya locus $\AA (R)$, so that $R/\mm R \cong  M_n (\KK)$. Since $\tr_{\mm}$ is
non-zero (by definition of an almost representation theoretic trace) and the kernel of the induced bilinear form is an ideal, the form
 $\lcor -, - \rcor_{\tr_{\mm}}$ is non-degenerate. This implies that $D_{n^2}(R/Z(R),\tr)$ is not contained in $\mm$.

Suppose on the other hand that $\mm$ is not in the Azumaya locus. Then, by PI theory,
(see \thref{char-Azum}(b), the key being Mueller's theorem [21]), $R/\mm R$ is not semisimple.
By hypothesis the trace $\tr$ is representation theoretic on $R/\mm R$, from which it follows,
essentially by linear algebra (\leref{2nd}), that \eqref{induce} holds, so that $\lcor -, - \rcor_{\tr}$
induces a symmetric bilinear form on $\overline{R} := (R/\mm R)/J(R/\mm R)$.
But now some further PI-theory, namely a version (\thref{char-Azum}(d)) of the
Additivity Principle due to Braun \cite{B}, implies that $\mathrm{dim}_{\KK} (\overline{R}) < n^2$. From this it follows
that $D_{n^2}(R/Z(R), \tr) \subseteq \mm$.

\subsection{ Historical remark}\label{history} The above sketch indicates that the underlying idea of the
proof is a very old one: namely, it is to use the nondegeneracy or otherwise of a symmetric
bilinear form induced by a representation to test for semisimplicity of a finite dimensional
algebra. Since the nondegeneracy of the form is determined by the determinant of the
discriminant, we arrive at the so-called \emph{discriminant test for semisimplicity} of a finite dimensional algebra $A$. When $A$ is commutative, the validity of this test was proved by
Weierstrass \cite{W} in 1884; the extension to noncommutative algebras was obtained by Molien \cite{Mo} in
1893. For further information on the relevant history, see Curtis, \cite[pages 54-55 and footnote, page 189]{C}.

\subsection{Notation}
\label{notate}
Throughout the paper $\KK$ will be an algebraically closed field, and $R$ will denote a $\KK$-algebra which is a finite module over its affine center $Z(R)$. Many of our
results could be formulated and proved for a more general field $\KK$, but for brevity and for the avoidance of technical complications we have not considered such
generalisations here. Given a maximal ideal $\mm$ of $Z(R)$, we denote the set of isomorphism classes of simple $R/\mm R$-modules by $\Irr_{\mm}(R)$.

For a commutative affine $\KK$-algebra $A$ and a subset $S \subseteq A$,
$\KK \lcor S \rcor $ will denote the unital subalgebra of $A$ generated by $S$, and $\lcor S \rcor$
will denote the ideal of $A$ generated by $S$. The group of units of $A$ will be denoted by $A^\times$. For an ideal $I$ of $A$, $\VV(I)$ will denote the zero locus of $I$ in $\maxspec A$. The Jacobson radical of an algebra $B$ is denoted $J(B)$.
Finally, set $\Zset_+ = \{1,2, \ldots \}$.
\medskip
\\
\noindent
{\bf Acknowledgements.} We are indebted to Amiram Braun who proposed an improvement of the first version of the Main Theorem
to all fields without a restriction on their characteristics.
We are thankful to Jason Bell, Ken Goodearl, Pavel Etingof, Tom Lenagan,
Jesse Levitt, Bach Nguyen, Manuel Reyes, Kurt Trampel, Daniel Sternheimer and James Zhang for very helpful discussions and comments on the first
version of the manuscript.
\sectionnew{Basic definitions}
\label{disc}
\subsection{Trace maps - definitions}
\label{2.1}

Concerning traces we shall make frequent use of slightly non-standard terminology, for which we need to consider first the case of a finite dimensional $\KK$-algebra $B$ and a finite dimensional $B$-module $V$. The composition
\[
\tr_V \colon B \to \End_\KK(V) \stackrel{\tr}\lra \KK,
\]
where the last map is the standard trace map on endomorphisms of finite dimensional vector spaces, is a trace-like map 
in the sense of \deref{trace}(1), and is a trace map whenever it is non-zero. Without loss of generality $V$ 
can be assumed to be semisimple; for, if $0 = V_0 \subset V_1 \subset \cdots \subset V_n = V$ is a composition series of $V$ 
and $\widehat{V}$ is defined to be $\oplus_{i=1}^n V_i/V_{i-1}$, then, clearly,
\[
\tr_V = \tr_{\widehat{V}}.
\]
\bde{trace}  Let $R$ be as in $\S$\ref{notate}, and let $C$ be a $\KK$-subalgebra of $Z(R)$ over which $R$ is a finite module.
\begin{enumerate}
\item A {\em{trace-like map}} from $R$ to $C$ is a map
$\tr \colon R \to C $ that has the properties:
\begin{enumerate}
\item[(i)] {\em{(}}$C$-linearity{\em{)}} $\tr(zx + wy) = z \tr(x) + w \tr(y)$ for $x,y \in R$, $z,w \in C$, and
\item[(ii)] {\em{(}}cyclicity{\em{)}} $\tr( xy) = \tr(yx)$ for $x, y \in R$.
\end{enumerate}
\item A \emph{trace} from $R$ to $C$ is a non-zero trace-like map from $R$ to $C$.

\item An \emph{almost representation theoretic trace-like map} [resp. {\em{almost representation theoretic trace}}] is a trace-like map 
[resp. trace] $\tr:R \longrightarrow Z(R)$ such that,  for all
$\mm \in \maxspec Z(R)$, there
exists a function $s_\mm : \Irr_\mm(R) \to \KK$  such that
the following diagram commutes:
\[
\begin{CD}
R             @>\tr>>                                Z(R)                                \\
@VVV                                                 @VVV                             \\
R/\mm R @>\sum_{V \in \Irr_{\mm}(R)}s_\mm (V) \tr_{V}>>    Z(R)/\mm \cong \KK.
\end{CD}
\]
\item A \emph{representation theoretic trace} is an almost representation theoretic trace $\tr:R \longrightarrow Z(R)$  with $s_{\mm}(V)$ non-zero for the unique $V \in \Irr_{\mm}(R)$, for all $\mm \in \mathcal{A}(R)$.
\end{enumerate}
\ede

To ensure that certain trace maps defined below take values in $Z(R)$, we will often need to impose the condition that $R$
coincides with its {\em{trace ring}} $T(R)$, whose definition we now recall. Assume that $R$ is prime, and let $Q$ be the fraction field of $Z(R)$. 
Thus $R \otimes _{Z(R)} Q$, the ring of fractions of $R$, is a central simple $Q$-algebra,
by (a special case of) Posner's theorem, \cite[$\S$I.13.3]{BG}.  Recall that $n$ denotes the PI degree of $R$; 
that is, $\dim_Q (R \otimes_{Z(R)} Q) = n^2$, \cite[$\S$I.13.3]{BG}.
Since $R$ is a finitely generated torsion-free $Z(R)$-module, it is, in the terminology of \cite[Chapter 2, $\S$8]{Re}, a
\emph{$Z(R)$-order} in $R \otimes _{Z(R)} Q$.  There exists a finite field extension $F$ of $Q$ which splits
$R \otimes _{Z(R)} Q$; namely,
$R \otimes_{Z(R)} F \cong M_n(F)$, see  \cite[\S 7b]{Re}.

Now, for $1 \leq i \leq n$, one considers the compositions
\begin{equation}\label{split}
c_i : R \hra R \otimes_{Z(R)} Q \hra R \otimes_{Z(R)} F \cong M_n(F) \stackrel{\sig_i}{\lra} F.
\end{equation}
Here $\sig_i : M_n(F) \to F$ are the elementary symmetric functions of the eigenvalues of a matrix:
\[
\sig_i(A) = \sum_{j_1 < \ldots < j_i} \la_{j_1} \ldots \la_{j_i}
\]
where $\la_1, \ldots, \la_n$ are the eigenvalues in $\ol{F}$ of $A \in M_n(F)$.
One easily shows that
$c_i : R \to Q$. The {\em{trace ring}} of $R$ is defined by
\begin{equation}
\label{trace-ring}
T(R) := R\langle c_i(r), r \in R, 1 \leq i \leq n \rangle \subset R \otimes_{Z(R)} Q.
\end{equation}
If $Z(R)$ is normal, then $R = T(R)$, see \cite[Theorem 10.1]{Re}.

Furthermore, one defines the characteristic polynomial of $a \in R$ by
\[
\chi_a(x) = x^n - c_1(a) x^{n-1} + \ldots + (-1)^n c_n(a) \in Q[x].
\]
We have 
\begin{equation}
\label{prop-red}
\chi_a(a) =0 \; \; \mbox{for all} \; \; a \in R.
\end{equation}

\subsection{Trace maps - examples} \label{2.1extra} 
Here are some canonical examples of trace maps. Throughout, $R$ satisfies the hypotheses of $\S$\ref{notate}.

\begin{enumerate}
\item Suppose that $R$ is a free $Z(R)$-module of finite rank $m$, and consider the composition
\[
R \; \hra \End_{Z(R)} (R) \stackrel{\tr}{\lra} Z(R)
\]
where the first map is given by left multiplication and the second is the usual trace of $m \times m$ matrices. Choosing a $Z(R)$-basis of $R$ yields an isomorphism of $\End_{Z(R)} (R)$ with $M_m (Z(R))$, and
then the last map is the usual trace of matrices. This composition is
called the {\em{regular trace}} of $R$, denoted by $\tr_{\reg}$. It is a trace map in the sense of the definition in $\S$\ref{2.1} provided $\charr \KK \nmid m$. It is then clearly representation theoretic, with $V_{\mm} = R/\mm R$ for all $\mm \in \mathrm{Maxspec}Z(R)$.

\bigskip

\item Suppose that $R$ is prime with PI-degree $n$ and that $R = T(R)$. Let $Q$ be the fraction field of $Z(R)$, and consider the following composition of maps, the case $i = 1$ of \eqref{split},
\begin{equation}
\label{tr-red}
 R \hra R \otimes_{Z(R)} Q \hra R \otimes_{Z(R)} F \cong M_n(F) \stackrel{\tr}{\lra} F,
\end{equation}
where the last map is the usual trace. It is called the \emph{reduced trace} and will be denoted $\tr_{\red}$. By hypothesis and construction the reduced trace $\tr_{\red} = c_1$ and indeed \emph{all} maps $c_i$ take values in $Z(R)$. It is \emph{not} obvious, but it is true and will be shown in \prref{re}, that
$$ \textit{ for all fields } \KK, \, \tr_{\red} \textit{ is a representation theoretic trace in the sense of } \S 2.1.$$

\bigskip

\item Let $R$, $Z(R)$, $Q$ and $n$ be as in (2), and consider the composition
\[
\tr_{\st} : R \; \hra \End_{Q} (R \otimes _{Z(R)} Q) \stackrel{\tr}{\lra} Q,
\]
where the first map is again given by left multiplication and the second is again the usual trace map.
This composition is called the {\em{standard trace}} of $R$. It is given in terms of the reduced trace by
\begin{equation}
\label{stan-red}
\tr_{\st} = n \tr_{\red}.
\end{equation}
Indeed, in terms of the splitting field $F$ in (2), $\tr_{\st}$ is given by the composition
\begin{equation}
\label{tr-stan2}
\tr_{\st} : R \hra R \otimes_{Z(R)} Q \hra R \otimes_{Z(R)} F \cong M_n(F) \stackrel{\tr_{\reg}}{\lra} F,
\end{equation}
where $\tr_{\reg} : M_n(F) \to F$ is the regular trace on $M_n(F)$. It is straightforward to check that
that the usual trace $\tr : M_n(F) \to F$ and the regular trace $\tr_{\reg} : M_n(F) \to F$ satisfy
$\tr_{\reg} = n \tr$. Thus \eqref{stan-red} follows from \eqref{tr-red} and \eqref{tr-stan2}.

If $R = T(R)$ (and in particular, if $Z(R)$ is normal), the standard trace takes values in $Z(R)$. In this case $\tr_{\st}$ is a trace in the sense of $\S$\ref{2.1} if $\charr \KK \nmid n$, and is representation theoretic when $\charr \KK \notin [1, n]$. 

\bre{standard} 
In view of \eqref{stan-red}, all facts that we prove for the reduced trace are
valid also for the standard trace under the additional assumption that $\charr \KK \notin [1, n]$.
\ere
\bigskip

\item Let $n$ be a positive integer and let $\KK$ be a field such that $\charr \KK \notin [1,n]$. Let $R$ be a $\KK$-algebra with trace $\tr : R \to C$, where $C$ is a central subalgebra.
Following \cite{P}, for $a \in R$, define its
$n$-th characteristic polynomial $\chi_{n,a}(x) \in Z(R)[x]$ to be
\[
\chi_{n,a}(x) := x^n - \wh{c}_1(a) x^{n-1} + \cdots + (-1)^n \wh{c}_n(a)
\]
where the functions $\wh{c}_i : R \to C$ are defined as follows. Write the $i$-th elementary
symmetric function $\sig_i$ in $n$ indeterminates $\la_1, \ldots, \la_n$ in terms of the Newton power sum functions
$\psi_j:=\la_1^j + \cdots + \la_n^j$ as
\[
\sig_i = p_i( \psi_1, \ldots, \psi_i) \quad \mbox{for} \quad p_i(t_1, \ldots, t_i) \in \Zset[(i!)^{-1}][t_1, \ldots, t_i].
\]
Then set $\wh{c}_i(a) := p_i( \tr(a), \ldots, \tr(a^i))$.

A Cayley--Hamilton algebra of degree $n$, as defined by Procesi \cite{P}, is an affine $\KK$-algebra $R$ with trace $\tr : R \to C$ such that
\[
\chi_{n,a}(a) =0 \; \; \mbox{for all}  \; \; a \in R \quad \mbox{and} \quad \tr(1) =n.
\]
We refer the reader to \cite[Sect. 4]{DP} and \cite[Sect. 1]{L} for detailed expositions.
In the special case when
\[
\charr \KK \notin [1,n]
\]
the reduced trace $\tr_{\red} : R \to Z(R)$ of (2) fits this framework because of \eqref{prop-red}.

\bigskip

\end{enumerate}


\subsection{Discriminants and discriminant ideals}
\label{2.1A}

Let $R$ be as stated in $\S$\ref{notate}, and let $C$ be a subalgebra of $Z(R)$ over which $R$ is a finite module. Let $\tr$ be a trace map from $R$ to $C$, and fix a
positive integer $m$. The {\em{$m$-discriminant ideal}} $D_m(R/C,\tr)$ of $R$ over $C$
with respect to $\tr$ is the ideal of $C$ generated by
\[
\det( [\tr(y_i y_j)]_{i,j=1}^m) \quad \mbox{for all $m$-tuples} \quad (y_1, \ldots, y_m) \in R^m.
\]
The {\em{modified $m$-discriminant ideal}} $MD_m(R/C,\tr)$ of $R$ is the ideal of $C$ generated by
\[
\det( [\tr(y_i y'_j)]_{i,j=1}^m) \quad \mbox{for all pairs of $m$-tuples} \quad (y_1, \ldots, y_m), (y'_1, \ldots, y'_m) \in R^m,
\]
see \cite[Definition 1.2 (2)]{CPWZ}. Clearly,
\begin{equation}
\label{incl}
D_m(R/C,\tr) \subseteq MD_m(R/C,\tr).
\end{equation}

Continuing with the above notation, let
\[
(\tilde{y}_1, \ldots, \tilde{y}_m), (\tilde{y}'_1, \ldots, \tilde{y}'_m) \in R^m \quad \mbox{and}  \quad (y_1, \ldots, y_l), (y'_1, \ldots, y'_l) \in R^l, \quad l \geq m,
\]
be such that
\[
(\tilde{y}_1, \ldots, \tilde{y}_m) = (y_1, \ldots, y_l) z \quad \mbox{and} \quad
(\tilde{y}'_1, \ldots, \tilde{y}'_m) = (y'_1, \ldots, y'_l) z'
\]
for some $z, z' \in M_{l \times m}(C)$. By row and column expansions one gets that
\begin{equation}
\label{change}
\det( [\tr(\tilde{y}_i \tilde{y}'_j)]_{i,j=1}^m) = \sum_{
\tiny{
\begin{matrix}
I, J \subseteq[1,l]
\\
| I | = | J |=m
\end{matrix}
}
}
\det (z_I) \det (z'_J) \det( [\tr(y_i y'_j)]_{i \in I, j \in J})
\end{equation}
where $z_I$ and $z'_J$ are the $I \times [1,m]$ and $J \times [1,m]$ minors of $z$ and $z'$. This implies part (a) of the following lemma, useful for computing
modified discriminant ideals.

\ble{gen} Let $R$, $C$, $\tr$ and $m$ be as above. Let $L$ be a subset of $R$ such that
\[
R = C L.
\]
\begin{enumerate}
\item[(a)] \cite[Lemma 1.11(1)]{CPWZ} The modified $m$-discriminant ideal $MD_m(R/C,\tr)$ is generated as an ideal of $C$ by
\[
\det( [\tr(y_i y'_j)]_{i, j =1}^m)
\]
for the $m$-element subsets $\{y_1, \ldots, y_m\}$ and $\{y'_1, \ldots, y'_m\}$ of $L$.
\item[(b)] Let $(w_1, \dots , w_m) \in R^{m}$ and suppose there are elements $z_1, \ldots , z_m$ of $C$, at least one of which is regular and such that
$\sum_i z_i w_i = 0.$ Then, for all $(u_1, \ldots , u_m) \in R^{m}$, $\det( [\tr(w_i u_j)]_{i, j =1}^m)= 0.$

\item[(c)] Suppose that every list of elements of $L$ of cardinality $m$ satisfies a non-trivial $C$-linear relation in the sense of part (b).
Then $MD_m(R/C,\tr) = \{0\}.$

\item[(d)] If $|L| <m$, then $MD_m(R/C,\tr) =0$.

\item[(e)] $MD_1(R/C,\tr)= \tr(R).$
\end{enumerate}
\ele

\begin{proof}(b) With the stated hypothesis on $\{w_i\}$, for each fixed $j = 1, \dots , m$,
\[
\sum_i z_i \tr (w_i u_j) = \tr((\sum_i z_i w_i)u_j) = 0,
\]
so, after an appropriate localization by a regular element of $C$, one of the rows of the matrix $[\tr(w_i u_j)]_{i, j =1}^m$
can be expressed in terms of the others. Therefore the determinant of this matrix equals 0.

(c) This follows from (a) and (b).

(d) is a special case of (c), and (e) is trivial.
\end{proof}

\bco{zero} Let $R$ be a prime $\KK$-algebra which is a finite module over its affine center $Z(R)$, let $\tr$ be a trace map from $R$ to $Z(R)$, and let $n$ be the PI
degree of $R$. Then $MD_m(R/Z(R),\tr) = 0$ for all $m > n^2$.
\eco
\begin{proof} Let $Q$ be the quotient field of $Z(R)$, so that $\mathrm{dim}_Q (R \otimes_{Z(R)} Q) = n^2$ by Posner's Theorem, \cite[\S I.13.3]{BG}. The corollary
follows from \leref{gen}(c) with $C = Z(R)$ and $L$ any finite generating set for $R$ as a $Z(R)$-module, since for $m > n^2$, every subset of $L$ of cardinality
$m$ is linearly dependent over $Q$.
\end{proof}

Recall that two elements $a, b \in A$ are \emph{associates}, denoted
\[
a =_{A^\times} b,
\]
if $a =u b$ for some $u \in A^\times$. Continuing with $R$, $C$ and $\tr$ as above, suppose that $R$ is a free $C$-module of finite rank $m$, with basis $\{y_1,, \ldots , y_m \}$. The \emph{discriminant}
(of $R$ over $C$ with respect to $\tr$) is defined to be the element
\[
d(R/C,\tr) :=_{C^\times} \det([\tr(y_iy_j)]_{i,j=1}^m) \in C.
\]
By a special case of \eqref{change}, this is independent of the choice of basis, meaning that for different
$C$-bases of $R$ the right hand sides above are associates of each other.

\ble{poly} Let $R$ be a $\KK$-algebra which is a finite module over its affine central subalgebra $C$, and let $\tr$ be a trace map from $R$ to $C$.
\begin{enumerate}
\item[(a)] \cite[Lemma 1.4(5)]{CPWZ} Let $m$ and $n$ be positive integers with $n \leq m$. Then
\[
MD_m(R/C,\tr) \subseteq MD_n(R/C,\tr).
\]
\item[(b)] Suppose that $R$ is a free $C$-module of rank $m$. Then
\begin{equation}
\label{princ}
MD_m(R/C,\tr)= D_m(R/C,\tr) = \langle d(R/C,\tr)\rangle.
\end{equation}
\item[(c)] Suppose that $C$ is a domain and $R$ is a projective $C$-module of rank $m$. Then the first equality in \eqref{princ} remains valid.
\item[(d)] If $\mathrm{(c)}$ holds, then $D_m(R/C, \tr)$ either equals $0$ or $C$, or else every prime ideal of $C$ minimal over $D_m(R/C, \tr)$ has height 1.

\end{enumerate}
\ele
\begin{proof} (b) Both equalities in \eqref{princ} are easy consequences of \leref{gen}(a). The second one is \cite[Theorem 10.2]{Re}.

(c) If the inclusion \eqref{incl} is strict, then it remains strict after localising at a suitable maximal ideal $\mm$ of $C$. But, denoting the extension of $\tr$ to $R_{\mm}$
by $\tr_{\mm}$, we find that $D_m(R_{\mm}/C_{\mm},\tr_{\mm}) = D_m(R/C,\tr)R_{\mm}$, by \cite[Exercise 4.13]{Re}, and an analogous result applies to 
$MD_m(R/C,\tr_{\mm})$. Since $R_{\mm}$ is $C_{\mm}$-free of rank $m$, this contradicts (a).

(d) Suppose that (c) holds and that $D_m(R/C,\tr)$ is a proper non-zero ideal. Let $\mathfrak{p}$ be a prime ideal of $C$ minimal over  $D_m(R/C,\tr)$, and localise at
$\mathfrak{p}$. Then $\mathfrak{p}$ has height at most one by \eqref{princ} and the Principal Ideal Theorem, \cite[Theorem 10.1]{E}.
\end{proof}
\bre{normal}
Suppose that $R$ satisfies the hypotheses of the Main Theorem, and is moreover a Cohen--Macaulay $Z(R)$-module. Take $C$ to be a polynomial subalgebra of
$Z(R)$ over which $R$ is a finite module, noting that such $C$ exist by Noether normalisation. Then $R$ is a free $C$-module,
by \cite[Theorem 3.7]{BrMac}, and the lemma applies.
\ere
\subsection{Examples of discriminant ideals}
\label{discex}
This subsection gives a number of examples of discriminant ideals which illustrate aspects of the above lemmas and which will also help to determine the limits of
possible extensions of results proved later in the paper.
\\

\noindent (1) This example shows that, notwithstanding \leref{poly}(b) and (d), neither the ideal $MD_m(R/C,\tr)$ nor $D_m(R/C,\tr)$, even when proper, need be contained in a
height one prime ideal.

Take $Z = \KK [x,y]$ with $\charr \KK \neq 2$, and let $R$ be the subalgebra of $M_2(Z)$ given by
\[
R :=
\begin{pmatrix}
Z &\langle x,y \rangle \\
Z &Z
\end{pmatrix}.
\]
Thus $Z(R) = Z$, $R$ has PI-degree 2, and
\[
R = Ze_{11} + Ze_{21} + Ze_{22} + Zxe_{12} + Zye_{12}.
\]
Calculating $\mathrm{det}([\tr_{\reg}(y_iy_j')]_{i,j=1}^5)$ using the above five $Z$-module generators for $\{y_i\}$ and $\{y_j'\}$ gives $0$. By \leref{poly}(b) or by
\coref{zero},
\[
MD_m(R/Z,\tr_{\reg}) = D_m(R/Z,\tr_{\reg}) =0 \quad \mbox{for} \quad m \geq 5.
\]
On the other hand, using in turn the two $4-$tuples of elements of $R$ got by omitting either $xe_{12}$ or $ye_{12}$ from the list of five generators gives
\[
MD_4(R/Z,\tr_{\reg}) = D_4(R/Z,\tr_{\reg}) = \langle x,y \rangle,
\]
in line with the prediction of the Main Theorem since the Azumaya locus of $R$ is easily checked to be
$\Aset^2 \backslash \{(0,0)\}$.
\\
\medskip

\noindent (2)\cite[Example 1.3(3)]{CPWZ} This example illustrates the following noteworthy features.
\begin{itemize}
\item The irreducible components of the algebraic set determined by the modified $m$-discriminant ideal need not be of equal dimension;
\item The ascending chain of discriminant ideals guaranteed by Lemmas \ref{lpoly}(a), \ref{lgen}(d) and \ref{lgen}(e) can have more than 3 distinct terms.
\end{itemize}
\medskip

Let $\KK$ have characteristic other than 2 or 3, and let $R$ be the quantum polynomial algebra $\KK_{\bf{p}}[X_1,X_2,X_3]$,
where $X_iX_j = p_{ij}X_jX_i$ and $p_{12}= p_{13}= -1, p_{23}= 1.$  In this case,
\[
\mathrm{gl.dim}(R) = \mathrm{Krull dim}(R) = \mathrm{GKdim}(R) = 3;
\]
here, the value for the global dimension follows from \cite[Theorem 7.5.3(iii)]{McR}, while the fact that $R$ is a finite module over the central polynomial algebra $\KK [ X_1^2, X_2^2,X_3^2 ]$ gives the remaining dimensions. By a direct calculation, or see for example \cite[Lemma 2.2(2)]{CPWZ},
\[
Z(R) =\KK \langle X_1^2, X_2^2,X_3^2, X_2X_3 \rangle.
\]
In particular, it is easy to check that $Z(R)$ has infinite global dimension. More precisely, the
singular locus $\mathcal{S}(Z(R))$ of $Z(R)$ is the subset
\begin{equation}
\label{sing}
\mathcal{S}(Z(R)) = \{\mm \in \mathrm{maxspec}(Z(R)) : \langle X_2^2, X_3^2 \rangle \subset \mm\}.  \end{equation}
Moreover, the argument used in the proof of \coref{gldim} below shows that $R$ is not a projective $Z(R)$-module. Setting $Q := Q(Z(R))$, a straightforward check shows that $R$ has PI-degree 2, equivalently that $\mathrm{dim}_Q(Q \otimes_{Z(R)} R) = 4$.
A minimal generating set for $R$ as a $Z(R)$-module is
\[
\{ 1, X_1, X_2, X_3, X_1X_2, X_1X_3 \}.
\]
From \cite[Example 1.3(3)]{CPWZ} we obtain that
\begin{equation}
\label{mixed}
MD_4(R/Z(R), \tr_{\st}) = \langle X_1^4X_2^iX_3^{4-i} : 0 \leq i \leq 4 \rangle,
\end{equation}
and a similar calculation shows that
\begin{equation}
\label{small}
MD_3(R/Z(R), \tr_{\st}) = \langle X_1^2X_2^2, X_1^2X_3^2,X_1^2X_2X_3  \rangle.
\end{equation}
From \eqref{mixed} and \eqref{small}, the algebraic subset of $\mathrm{maxspec}(Z(R))$ defined by both ideals $MD_4(R/Z(R), \tr_{\st})$ and $MD_3(R/Z(R), \tr_{\st})$ is
\begin{equation}
\label{unpure}
\mathcal{V}(\langle X_1^2 \rangle \cap \langle X_2^2, X_3^2 \rangle).
\end{equation}
By the Main Theorem, this algebraic set constitutes the non-Azumaya locus of $R$. The first bullet point above follows from \eqref{unpure}.
The second bullet point is immediate from \eqref{mixed} and \eqref{small}.
\sectionnew{Proof of the main result}
\label{pf}
\subsection{Preliminaries on primes}\label{3.1} The proof of part (a) of the Main Theorem relies on the characterization of the Azumaya locus of a PI algebra given in part (b)(vi)
of \thref{char-Azum} below. For its proof we need to recall that prime ideals $P$ and $Q$ of a noetherian ring $T$ are said to be
({\em{second layer}}) {\em{linked}} if either
 $P \cap Q/PQ$ is not a torsion $R/P-R/Q-$bimodule, or $P\cap Q/QP$ is not a torsion $R/Q-R/P-$bimodule. Letting $\sim$ denote the equivalence relation on
 $\mathrm{Spec}(T)$ generated by the second layer links, the $\sim-$classes are called the \emph{cliques} of prime ideals of $T$, generalising the blocks of simple
modules of a finite dimensional algebra. For further details and references, see \cite[Chapters 12-14]{GW}.

Continue throughout the remainder of $\S$3 to assume that the following hypothesis (H) holds:
\begin{align*}\label{hypo}\textbf{(H) }&R \textit{ is a prime affine }\KK-\textit{algebra which is a}\\ &\textit{finite module over its noetherian center }Z(R).
\end{align*}
By the Artin--Tate lemma, \cite[Lemma 13.9.10]{McR},  $Z(R)$ is an affine $\KK$-algebra, so that $Z(R)$ and hence $R$ are noetherian rings. By Kaplansky's
theorem \cite[I.13.3]{BG} every primitive ideal $M$ of $R$ is maximal, so that $R/M \cong M_t(\KK)$ for some $t \geq 1$, and $M \cap Z(R)$ is a maximal ideal of
$Z(R)$. Conversely, given a maximal ideal $\mm$ of $Z(R)$, there is a finite but non-empty set of maximal ideals $M_1, \dots , M_s$ of $R$
for which $M_i \cap Z(R) = \mm$.

\subsection{Characterisations of the Azumaya locus}\label{azum}  Recall that hypothesis (H) from $\S$\ref{3.1} is in force. A prime ideal $P$ of $R$ is called
\emph{regular} if the PI-degree of $R/P$ equals that of $R$ (in general the former
is less than or equal to the latter). A ring $A$ is called an \emph{Azumaya algebra} over its center $Z(A)$, \cite[13.7.6]{McR}, \cite[III.1.3]{BG}, if $A$ satisfies the two
conditions
\begin{itemize}
\item $A$ is a finitely generated projective $Z(A)$-module; and
\item the ring homomorphism
\[
A \otimes_{Z(A)} A \to \End_{Z(A)} (A) \quad \mbox{given by} \quad
a \otimes b \mt (x \mt axb)
\]
is an isomorphism.
\end{itemize}
The \emph{Azumaya locus} $\AA(R)$ of a PI algebra $R$, satisfying hypothesis (H) of $\S$\ref{3.1}, is defined
to be the set of $\mm \in \maxspec Z(R)$ satisfying any of the equivalent conditions (b)(i)--(vi) of the following theorem. Note that $\AA(R)$ is a proper closed subset of
$\mathrm{Maxspec}Z(R)$ by \cite[Theorem III.1.7]{BG}. Recall that the additional hypothesis $R = T(R)$ imposed in (c)-(f) of the following result holds, for example, when $Z(R)$ is normal \cite[Theorem 10.1]{Re}.

\bth{char-Azum} Assume hypothesis {\em{(H)}}, that $R$ is a prime affine $\KK$-algebra which is
a finitely generated module over its center $Z(R)$. Assume that $\KK$ is algebraically closed. Let $n$ be the PI degree of $R$. Let $M$ be a maximal ideal of $R$, and
denote by $\mm:= M \cap Z(R)$ the corresponding maximal ideal of $Z(R)$.
\begin{enumerate}
\item[(a)] Let $\n$ be a maximal ideal of $Z(R)$. The set of maximal ideals $N_i$ of $R$ for which $N_i \cap Z(R) = \n$ has cardinality $t$, where $1 \leq t \leq n$.
\item[(b)] The following are equivalent:
\begin{enumerate}
\item[(i)] $M$ is a regular maximal ideal of $R$.
\item[(ii)] $R_\mm$ is an Azumaya algebra over $Z(R)_\mm$.
\item[(iii)] $M = \mm R$.
\item[(iv)] The unique irreducible $R/M$-module has maximal dimension amongst irreducible $R$-modules, namely $n$.
\item[(v)] $R/M \cong M_n(\KK)$.
\item[(vi)] $R/\mm R$ is a semisimple algebra.

\noindent Moreover, these equivalent conditions imply that
\item[(vii)] $\dim_{\KK}( (R/\mm R)/J(R/ \mm R)) = n^2$.
\end{enumerate}
\item[(c)] Suppose in addition that $R$ coincides with its trace ring $T(R)$. Then, in $(\mathrm{b})$, $(\mathrm{vii})$ is equivalent to $(\mathrm{i})-(\mathrm{vi})$.
\item[(d)] Continue to assume that $R=T(R)$, and suppose that $\mm$ is a maximal ideal of $Z(R)$ for which the equivalent conditions
$(\mathrm{i})-(\mathrm{vii})$ of $(\mathrm{b})$ do not hold. Let $t$ be the number of irreducible $R/\mm R$-modules, and let
their $\KK$-dimensions be $m_1, \ldots , m_t$. Then
\begin{equation}
\label{pence}
\dim_{\KK}( (R/\mm R)/J(R/ \mm R)) \; \leq \; n^2 - \mathrm{max}\{\sum_{i=1}^t m_i, 2 \} .
\end{equation}

\item[(e)] If $R = T(R)$ then $\dim_{\KK}((R/\mm R)/J(R/ \mm R)) \leq n^2$, and equality holds if and only if the equivalent conditions in $(\mathrm{b})$ are
satisfied.
\item[(f)] With the notation and hypotheses of $(\mathrm{d})$, the upper bound in $(\mathrm{a})$ can be strengthened to
\begin{equation} n \; \geq \; \sum_{i=1}^t m_i.
\end{equation}
\end{enumerate}
\eth
\begin{proof} (a) The lower bound is a consequence (after localising at $\mathfrak{n}$) of Nakayama's lemma; the upper bound is proved at \cite[Proposition 5]{B},
where the hypothesis imposed there, that $R = T(R)$, is not needed when the ground field is algebraically closed, (as assumed here).

(b) The equivalence of conditions (i)--(v) is proved in \cite[Theorem III.1.6]{BG}, and it is trivial that (iii) $\Rightarrow$ (vi) and that $(ii)\cup (v)$ $\Rightarrow$ (vii).

It remains to show that (vi) $\Rightarrow$ (iii). Let $\mm$ be a maximal ideal of $Z(R)$ such that $R/\mm R$ is semisimple, and let $M$ and $N$ be maximal ideals of
$R$ with $M\cap Z(R) = N\cap Z(R) = \mm$. By Muller's theorem, \cite{Mu}, \cite[Theorem III.9.2]{BG}, $M$ and $N$ belong to the same clique of $\maxspec(R)$.
However, thanks to hypothesis (vi), each of $M/ \mm R$ and $N/\mm R$ is generated by a central idempotent of $R/\mm R$. Hence, $M$ and $N$ are polycentral
ideals of $R$, so that $\mathrm{clique}(M) = \{M\}$ and $\mathrm{clique}(N) = \{N\}$ by \cite[Proposition 13.6]{GW} combined with \cite[Theorem 4.2.7]{McR}. Therefore
$M =N$, proving (iii).

(d), (f) Suppose that $R = T(R)$ and that $R_{\mm}$ is not an Azumaya algebra, so that (b)(i)-(vi) all fail for $\mm$, by (b). By \cite[Proposition 4]{B}, there are positive integers $z_i$, $(1 \leq i \leq t)$, such that
\begin{equation}
\label{ami}
n = \sum_{i=1}^t z_i m_i.
\end{equation}
In particular, $n \geq \sum_{i=1}^t m_i$, proving (f). It follows that
\begin{align}
\label{trump}
n^2 &\geq n\big(\sum_{i=1}^t m_i\big)\\
\nn
&= \sum_{i=1}^tm_i^2 + \sum_{i=1}^t (n - m_i)m_i \\
\nn
&\geq \sum_{i=1}^tm_i^2 + \sum_{i=1}^t m_i,
\end{align}
where the final inequality holds because $n - m_i > 0$ for all $i$, since hypothesis (b)(v) does not hold for any maximal ideal containing $\mm R$. Thus
\begin{equation}
\label{ivanka}
\dim_{\KK} (R/\mm R)/J(R/ \mm R) \; = \;   \sum_{i=1}^tm_i^2 \; \leq \; n^2 - \sum_{i=1}^t m_i.
\end{equation}
To complete the proof of \eqref{pence}, note finally that if $\sum_{i=1}^t m_i = 1$ and equality holds throughout \eqref{ivanka}, then $t = 1$ and $m_1 = 1$, so that
$\dim_{\KK}(R/\mm R)/J(R/ \mm R) = 1$ and \eqref{ivanka} becomes $1 = n^2 - 1$, which is impossible for $n \in \mathbb{Z}$.

(c) Suppose $R = T(R)$, and assume that (vii) holds for $\mm$. If $R_{\mm}$ is not Azumaya, then \eqref{ivanka} for $\mm$ contradicts (vii) for $\mm$. So $R_{\mm}$ is Azumaya.

Finally, (e) follows from (b), (c) and (d).
\end{proof}
\subsection{Trace forms for finite dimensional algebras}
\label{3.2}
Let $B$ be a finite dimensional algebra over a field $\KK$ and $\tr \colon B \to \KK$
be a trace map. Denote by  $\lcor -, - \rcor_{\tr} \colon B \times B \to \KK$ the  corresponding (symmetric,
bilinear) trace form, given by
\[
\lcor a, b \rcor_{\tr} := \tr(ab) \quad \mbox{for} \quad a, b \in B.
\]

Define the \emph{kernel} of the trace form   $\lcor -, - \rcor_{\tr}$ to be $L_{\tr} = \{b \in B : \tr (bx) = 0, \,  \forall \; x \in B \}$. It is clear that $L_{\tr}$ is an ideal of $B$; by definition,  $\lcor -, - \rcor_{\tr}$ is \emph{nondegenerate} if and only if $L_{\tr} = \{0\}$. Below and in the rest of $\S 3$ we use the notation $\tr_V$ introduced in 
\S\ref{2.1}.

\ble{tr-form} Fix notation as above.
\begin{enumerate}
\item[(a)] If $B \cong M_n(\KK)$ and $\tr \colon B \to \KK$ is any trace map, then
the trace form $\lcor -, - \rcor_{\tr}$ is nondegenerate.

\item[(b)] If $V_1, \ldots, V_j$ are non-isomorphic simple $B$-modules and $\tr \colon B \to \KK$ is given by 
\begin{equation}
\label{tr-assume}
\tr = s_1 \tr_{V_1} + \cdots + s_j \tr_{V_j}
\end{equation}
for some $s_j \in \KK^*$, then the kernel $L$ of the trace form $\lcor -, - \rcor_{\tr}$ equals
\[
L := \cap \mathrm{Ann}_B(V_i).
\]
In particular, $J(B) \subseteq L$.
\end{enumerate}
\ele

\begin{proof} (a) The kernel $L_{\tr}$ of the trace form is a proper two-sided ideal of $M_n(\KK)$ since $1 \notin L_{\tr}$.
The only such ideal of $M_n(\KK)$ is the zero ideal.

(b) It follows from assumption \eqref{tr-assume} that the trace map $\tr$ and the form $\lcor -, - \rcor_{\tr}$ descend to the algebra 
$B/J(B)$. Furthermore, \eqref{tr-assume} implies that $\tr$ is nonzero on the direct summands of the semisimple algebra $B/J(B)$ 
corresponding to the representations $V_1, \ldots, V_j$ and vanishes on the other direct summands. Part (b) follows from this property.
\end{proof}

\bre{degen} The conclusion of \leref{tr-form} (b) cannot be extended to trace forms which do not arise from a representation.
Indeed, consider any commutative, non-semisimple algebra $B$. Every non-zero functional $\tr \colon B \to \KK$ is a trace
map. This easily produces examples for which the corresponding trace form $\lcor -, - \rcor_{\tr}$ is non-degenerate.
\ere
\subsection{Descent of traces to semisimple algebras}
\label{main-sect}
Before proving the Main Theorem we need the following key lemma. The terminology concerning trace maps used here was introduced in $\S$\ref{2.1}.

Recall that for $\mm \in \maxspec Z(R)$, $\Irr_\mm(R)$ denotes the equivalence classes 
of the finite dimensional irreducible representations of $R$ which are annihilated by $\mm$. 
Recall from \eqref{trnote} that a trace map $\tr : R \longrightarrow Z(R)$ descends to a 
trace-like map (which might possibly be $0$),
\[
\tr_{\mm} \colon R/\mm R \to Z(R)/\mm \cong \KK.
\]
In terms of the canonical homomorphism
\begin{equation}
\label{chi}
\phi_\mm : R \to R/ \mm R,
\end{equation}
we have, for $a \in R$, 
\begin{equation}\label{bling} \tr_{\mm}(\phi_{\mm}(a)) = \phi_{\mm}( \tr (a)).
\end{equation}

\ble{2nd} Assume that $R$ is an affine $\KK$-algebra which is a finitely generated module over its center $Z(R)$. Let $\tr \colon R \to Z(R)$ be an almost 
representation theoretic trace-like map.
If $a \in R$ and $\mm \in \maxspec Z(R)$ are such that $a + \mm R \in J(R/\mm R)$, then
\[
\tr (a) \in \mm.
\]
In other words, $\tr_\mm$ descends to a trace-like map 
\begin{equation}
\label{desc}
\ol{\tr}_{\mm} \colon (R/\mm R) / J (R/\mm R) \to Z(R)/\mm \cong \KK
\end{equation}
on the semisimple algebra $(R/\mm R) / J (R/\mm R)$.
\ele

\begin{proof}
Since $\tr$ is an almost representation theoretic trace-like map, we can fix $s_\mm : \Irr_\mm(R) \to \KK$ satisfying \eqref{rep-theor-tr2}.
By \leref{tr-form}(b), \eqref{bling} and the definition of almost representation theoretic trace-like map,
\[
\phi_\mm (\tr(a)) = \tr_{\mm}(\phi_{\mm} (a)) =\sum_{ V \in \Irr_\mm(R)} s_\mm(V) \tr_V(\phi_ \mm (a))  = 0.
\]
Therefore $\tr (a) \in \mm R \cap Z(R) = \mm$, as required.
\end{proof}
\subsection{The reduced trace is representation theoretic}\label{reduce} In this subsection we show both that the reduced trace \emph{is} 
a trace in the sense of $\S$\ref{2.1} - that is, it is non-zero; \emph{and} that it is representation theoretic. For part (a), the argument is a minor reinterpretation 
of \cite[Proposition 4]{B} of Braun. 

\bpr{re} Let $R$ be a prime affine algebra over an algebraically closed field $\KK$, which is
a finitely generated module over its center $Z(R)$. Assume that $R$ coincides with its trace ring $T(R)$.
Let $n$ be the PI degree of $R$.

\begin{enumerate}
\item[(a)] For every $\mm \in \maxspec Z(R)$ there exists a function
\[
k_\mm : \Irr_\mm(R) \to \Zset_+
\]
such that
\begin{equation}
\label{cond1}
n = \sum_{V \in \Irr_\mm(R)} k_\mm(V) \dim_\KK V.
\end{equation}
Moreover, writing $\pi$ for the canonical homomorphism from $\Zset$ to $\KK$,
\begin{equation}
\label{cond2}
(\tr_{\red})_\mm = \sum_{V \in \Irr_\mm(R)}\pi( k_\mm(V)) \tr_V
\end{equation}
where $\tr_V : \End_\KK(V) \to \KK$ denotes the usual trace.

\item[(b)] The reduced trace $\tr_{\red}$ is a representation theoretic trace.  
\item[(c)] The standard trace $\tr_{\st}$ is a representation theoretic trace provided $\charr \KK \notin [1,n]$.
\end{enumerate}
\epr

Note that \eqref{cond2} implies \eqref{cond1} if $\charr \KK \notin [1,n]$, since in that case $\tr_{\red}(1) = n$ and $\tr_V (1) = \dim_{\KK}V$ for all $V \in \Irr_{\mm}(R)$. But over an arbitrary field there is no direct implication between the two.

\begin{proof} (a) Fix $\mm \in \maxspec Z(R)$ and recall the notation introduced at \eqref{split} in $\S$\ref{2.1} and the canonical epimorphism $\phi_{\mm}$ of \eqref{chi}. We show first that the functions $c_i : R \to Z(R)$ descend to functions
$c_{i, \mm} : R/ \mm R \to Z(R)/ \mm \cong \KK$, which then extend to functions
\[
c_{i, \mm} : (R/  \mm R) \otimes_\KK A \to A
\]
for every commutative $\KK$-algebra $A$. Choose $v_1, \ldots , v_g \in R$ such that $\{\phi(v_1), \dots , \phi (v_g)\}$ is a $\KK$-basis of $R/\mm R$. In \cite[equation (1), p.434]{B}, Braun restates Amitsur's key formula  \cite[(4.1)]{A} as follows. Let $x_1, \ldots , x_g$ be commuting indeterminates. Define an equivalence relation on the monoid of all words in $\{v_1, \ldots , v_g \}$ by declaring words $p$ and $q$ to be equivalent if one is a cyclic permutation of the other, and say that a word $p$ is \emph{indecomposable} if it is not a power of a word of strictly smaller length. Let $\mathcal{S}$ denote the set of equivalence classes of indecomposable words of length at least 1 and at most $n$. Then, for $i = 1, \ldots, n$,
\begin{equation}\label{Shimson} c_i(\sum_{j=1}^g x_j v_j) = \sum_{\ell} h_{\ell i}(x_1, \ldots , x_g)k_{\ell i} (v_1, \ldots , v_g), 
\end{equation}
where the $h_{\ell i}$ are polynomials in $\{x_1, \ldots , x_g \}$ with coefficients $\pm 1$ and no constant terms, and the $k_{\ell i}$ are polynomials with coefficients $\pm 1$ in $\{c_j (p_t) : 1 \leq j \leq n,\,  p_t \in \mathcal{S}\}$.

A straightforward calculation using \eqref{Shimson} shows that, if $r_1, r_2 \in R$ and $z \in \mm$, then, for some $w \in Z(R)$,
\[
c_i (r_1 + z r_2) = c_i (r_1) + zw. 
\]
In particular,
\[
c_i (r_1 + z r_2) -  c_i (r_1) \in \mm, 
\]
and our claim regarding the descent of the functions $c_i$ follows.

In \cite[Proposition 4]{B} and its proof, specifically, equation (9),  Braun proves the existence of a function
$k_\mm : \Irr_\mm(R) \to \Zset_+$, satisfying \eqref{cond1} and having the property
\begin{equation}
\label{cond3}
c_{n, \mm}(b)= \prod_{V \in \Irr_\mm(R)} \det(\pi_V(b))^{k_\mm(V)} \quad \mbox{for all} \quad
b \in R/( \mm R) \otimes_\KK A
\end{equation}
where the representation $\pi_V : R/ \mm R \to \End_\KK(V)$ is extended to a homomorphism
\[
\pi_V : (R/ \mm R) \otimes_\KK A \to \End_A(A \otimes_\KK V).
\]
For $a \in R/(\mm R)$, we have, using \eqref{cond1},
\[
c_{n, \mm}(x - a) = x^n - (\tr_{\red})_\mm(a) x^{n-1} + \; \mbox{lower order terms}.
\]
Moreover,
\begin{align*}
\prod_{V \in \Irr_\mm(R)} \det(\pi_V(x-a))^{k_\mm(V)} = x^n - \bigg{(} \sum_{V \in \Irr_\mm(R)} k_\mm(V) \tr_V(\pi_V(a) )&\bigg{)} x^{n-1}\\
&+ \; \mbox{lower order terms}.
\end{align*}
Apply \eqref{cond3} for $A = \KK[x]$ where $x$ is an indeterminate to identify the above two expressions. Thus \eqref{cond2} now follows from \eqref{cond3}.

(b) It is clear from \eqref{cond2} that $\tr_{\red}$ is a representation theoretic trace-like map. Let $\mm$ be in the Azumaya locus $\mathcal{A}(R)$. Since there exists a unique simple $R/\mm R$-module $V_\mm$, and $\dim_\KK V_\mm = n$, \eqref{cond1} shows that 
$$ k_\mm (V_\mm) = 1. $$
Hence, by \eqref{cond2},
$$ (\tr_{\red})_\mm = \tr_{V_\mm},$$
proving (b).


(c) This is clear from the definition of $\tr_{\st}$. It also follows from (b) using \eqref{stan-red} and \reref{standard}.
\end{proof}

\subsection{Proof of the Main Theorem - I}
\label{parts-ac} Observe first that parts (b) and (c) of the Main Theorem are restatements of \prref{re}(b) and (c), and part (f) of the Main Theorem 
follows from \coref{zero} and \eqref{incl}. 
We shall deal with parts (d) and (e) in \S\ref{partd}. \prref{main} below proves part (a).

\bpr{main} Let $R$ be a prime affine algebra over the algebraically closed field $\KK$, which is a finitely generated module over its center $Z(R)$.
Suppose that $R = T(R)$ and let $n$ be the PI degree of $R$. Let $\tr : R \to Z(R)$ be a representation theoretic  trace map. Then
\[
\VV(MD_{n^2}(R/Z(R),\tr) ) = \VV(D_{n^2}(R/Z(R),\tr)) = \maxspec Z(R) \backslash \AA(R).
\]
\epr

\begin{proof} First we show that
\begin{equation}
\label{incl2}
\VV(MD_{n^2}(R/Z(R),\tr) ) \subseteq \VV(D_{n^2}(R/Z(R),\tr)) \subseteq \maxspec Z(R) \backslash \AA(R).
\end{equation}
The first inclusion is a special case of \eqref{incl}.
For the second inclusion, choose $\mm \in \AA(R)$, so, using the notation \eqref{chi}, $R/\mm R \cong M_n(\KK)$. Choose $y_1, \ldots, y_{n^2} \in R$ such that
$\phi_\mm(y_1), \ldots, \phi_\mm(y_{n^2})$ is a basis of $R/\mm R$. By \leref{tr-form}(a) and the hypothesis that $\tr$ is almost representation theoretic, the trace form
$\lcor-, -\rcor_{\tr_\mm}$ is nondegenerate. Thus,
\[
\phi_\mm \big( \det \big( [ \tr (y_i y_j)]_{i,j=1}^{n^2} \big) \big) =
\det( [ \lcor \phi_\mm(y_i), \phi_\mm(y_j) \rcor_{\tr_\mm}]_{i,j=1}^{n^2}) \neq 0,
\]
which implies that $\mm \not \in \VV(D_{n^2}(R/Z(R), \tr) )$. This proves that
\[
\VV(D_{n^2}(R/Z(R), \tr)) \cap \AA(R) = \varnothing,
\]
proving \eqref{incl2}.

It remains to prove that
\begin{equation}
\label{plex}
\VV(MD_{n^2}(R/Z(R),\tr) ) \supseteq \maxspec Z(R) \backslash \AA(R).
\end{equation}
Let $\mm \in \maxspec Z(R) \backslash \AA(R)$. By \leref{2nd}, there is an induced trace-like map $\overline{\tr}_{\mm}$ 
on $R/\mm R/J(R/\mm R)$, and hence an associated bilnear form $\lcor -, - \rcor_{\ol{\tr}_\mm}$ on $(R/\mm R)/J(R/\mm R)$. 
The bilinear form could be 0; but, in any case, by  \thref{char-Azum}(c) and (d),
\[
\dim_\KK \big( (R/\mm R)/J(R/\mm R) \big) < n^2.
\]
Therefore, for all $(y_1, \ldots, y_{n^2})$ and $(y'_1, \ldots, y'_{n^2}) \in R^{n^2}$,
\[
\phi_\mm \big( \det \big( [ \tr (y_i y'_j)]_{i,j=1}^{n^2} \big) \big) =
\det( [ \lcor \phi_\mm(y_i) \phi_\mm(y'_j) \rcor_{\ol{\tr}_\mm}]_{i,j=1}^{n^2}) =0.
\]
Hence, $\mm \in \VV(MD_{\ell}(R/Z(R), \tr))$, which proves \eqref{plex}.
\end{proof}
Taking into account \leref{tr-form}(a), one easily sees that the proof of the proposition establishes the following more general fact than 
part (a) of the Main Theorem, which is of independent interest.
\bpr{MTa} Let $R$ be a prime affine algebra over the algebraically closed field $\KK$, which is a finitely generated module over its center $Z(R)$.
Suppose that $R = T(R)$ and let $n$ be the PI degree of $R$. Let $\tr : R \to Z(R)$ be a trace that
descends to trace-like maps on the semisimple quotients $\ol{\tr}_\mm : (R/\mm R)/J(R/\mm R) \to \KK$ for all $\mm \in \maxspec Z(R)$. Then
\begin{multline*}
\VV(MD_{n^2}(R/Z(R),\tr) ) = \VV(D_{n^2}(R/Z(R),\tr)) =  
\\
= \left( \maxspec Z(R) \backslash \AA(R) \right) \cup 
\{ \mm \in \AA(R) \mid \ol{\tr}_\mm =0 \}.
\end{multline*}
\epr
By \leref{2nd} every almost representation theoretic trace map satisfies the assumptions if the proposition. 
\bex{MTa2} The statement of part (a) of the Main Theorem does not hold in the more general situation of almost representation 
theoretic traces. Indeed, let $\tr : R \to Z(R)$ be a representation theoretic trace and $z \in Z(R)$ be a non-zero element. Then 
\[
\tr' : R \to Z(R) \quad \mbox{given by} \quad \tr'(a) := z \tr(a), \; \; \forall \; a \in R
\]
is an almost representation theoretic trace which is not representation theoretic if 
\[
\VV(z) \cap \AA(R) \neq \varnothing.
\]
Applying \prref{MTa} to $\tr'$, gives that
\[
\VV(MD_{n^2}(R/Z(R),\tr') ) = \VV(D_{n^2}(R/Z(R),\tr')) = \big(\maxspec Z(R) \backslash \AA(R) \big) \cup \VV(z).
\]
\eex
\subsection{Proof of the Main Theorem - II}
\label{partd} 
\begin{proof}[Proof of part (d) of the Main Theorem]  Let $\tr : R \to Z(R)$ be an almost representation theoretic trace 
with associated functions $s_\mm : \Irr_\mm(V) \to \KK$. 
Fix an element $\mm$ of $\maxspec Z(R)$. Set
\[
d_\mm := \sum_{V \in \Irr_\mm(R), s_\mm(V) \neq 0} (\dim_\KK V)^2,
\]
so $d_\mm$ is a positive integer by definition of $\tr$. Let 
\[
\overline{(\cdot)}:R/\mm R \longrightarrow (R/\mm R)/J(R/\mm R) = \overline{R/\mm R}.
\]
denote the canonical epimorphism.
By \leref{2nd}, the induced trace $\tr_{\mm} : R/ \mm R \to \KK$ descends to a trace
\[
{\ol{\tr}}_\mm :( R/ \mm R) / J( R/ \mm R) = \ol{R/\mm R} \to \KK.
\]
By Artin-Wedderburn,
\[
(R/ \mm R) / J( R/ \mm R) \cong \prod_{V \in \Irr_\mm(R)} \End_\KK(V),
\]
and, in terms of the right hand side, ${\ol{\tr}}_\mm$ is given by
\[
\ol{\tr}_\mm( \{ a_V : V \in \Irr_\mm(R) \} ) = \sum_{V \in \Irr_\mm(R)} s_\mm(V) \tr_V(a_V) \quad
\mbox{for} \quad a_V \in \End_\KK(V).
\]
Therefore the kernel of $\lcor -,  - \rcor_{{\ol{\tr}}_\mm}$ is 
\[
\Ker \lcor -,  - \rcor_{\ol{\tr}_\mm} = \prod_{V \in \Irr_\mm(R), s_\mm(V)=0} \End_\KK(V)
\]
and the restriction of $\lcor -,  - \rcor_{\ol{\tr}_\mm}$ to 
\[
\big( (R/ \mm R) / J( R/ \mm R) \big) / \Ker \lcor -,  - \rcor_{({\ol{\tr}_{\red}})_\mm} \cong \prod_{V \in \Irr_\mm(R), s_\mm(V) \neq 0} \End_\KK(V)
\]
is non-degenerate. The dimension of the above algebra equals $d_\mm$. 
Let $\ell$ be a positive integer. Recall from  \eqref{bling} the canonical homomorhism $\phi_\mm : R \to R/\mm R$. If 
\[
\ell \leq d_\mm, 
\]
then there exist $(y_1, \ldots, y_{\ell})$ and $(y'_1, \ldots, y'_{\ell}) \in R^{\ell}$ such that
\[
\phi_\mm \big( \det \big( [ \tr (y_i y'_j)]_{i,j=1}^{\ell} \big) \big) =
\det( [ \lcor \ol{\phi_\mm(y_i)}, \ol{\phi_\mm(y'_j)} \rcor_{\ol{\tr}_\mm}]_{i,j=1}^{\ell}) \neq 0.
\]
If 
\[
\ell > d_\mm,
\]
then for all $(y_1, \ldots, y_{\ell})$ and $(y'_1, \ldots, y'_{\ell}) \in R^{\ell}$ 
\[
\phi_\mm \big( \det \big( [ \tr (y_i y'_j)]_{i,j=1}^{\ell} \big) \big) =
\det( [ \lcor \ol{\phi_\mm(y_i)}, \ol{\phi_\mm(y'_j)} \rcor_{\ol{\tr}_\mm}]_{i,j=1}^{\ell}) = 0.
\]
Therefore
\[
\VV(D_{\ell}(R/Z(R),\tr) )  = \VV(MD_{\ell}(R/Z(R),\tr))
\]
and a maximal ideal $\mm \in \maxspec Z(R)$ belongs to this zero locus if and only if 
\[
\ell > d_\mm.
\] 
\end{proof}
The following theorem is a stronger form of part (e) of the Main Theorem without any assumptions on the characteristic of the base field $\KK$. 

\bth{red-tr-main} Let $R$ be a prime affine algebra over an algebraically closed field $\KK$, which is
a finitely generated module over its center $Z(R)$. Assume that $R$ coincides with its trace ring $T(R)$.

For all positive integers $\ell$, the zero sets of the $\ell$-discriminant ideal and the modified $\ell$-discriminant ideal of $R$ coincide and are given in
terms of the irreducible representations of $R$ by
\begin{multline*}
\VV(D_{\ell}(R/Z(R),\tr_{\red}) )  = \VV(MD_{\ell}(R/Z(R),\tr_{\red}))
\\
=  \Big\{ \mm \in \maxspec Z(R) \mid
\sum_{V \in \Irr_\mm(R), \charr \KK  \nmid k_\mm(V)} (\dim_\KK V)^2  < \ell \Big\},
\end{multline*}
where the functions $k_\mm : \Irr_\mm(R) \to \Zset_+$ for $\mm \in \maxspec Z(R)$ 
are the ones from \prref{re}.
\eth
To obtain \thref{red-tr-main} we apply \prref{re}(b) that the reduced trace is representation theoretic and part (d) of the Main Theorem. 
The function $k_\mm : \Irr_\mm(R) \to \Zset_+$ from \prref{re}(a) gives rise to a function $s_\mm : \Irr_\mm(R) \to \KK$ satisfying 
\eqref{rep-theor-tr2} by using the natural ring homomorphism $\Zset \to \KK$. Therefore, for $\mm \in \maxspec(R)$ and $V \in \Irr_\mm(R)$,  
\[
s_\mm(V) = \pi(k_\mm(V)) \neq 0 \quad \mbox{if and only} \quad \charr \KK  \nmid k_\mm(V).
\]
Recall that $\pi : \Zset \to \KK$ denotes the canonical homomorphism. 
\subsection{Supplementary remarks}\label{Poisson-ord}
(a) \textbf{Poisson orders:} Let $R$ be a $\KK$-algebra which is a finitely
generated module over a central subalgebra $C$.
Denote by $\Der_\KK(R)$ the Lie algebra of $\KK$-algebra derivations of $R$.
The pair $(R,C)$ is called a Poisson order \cite[Definition 2.1]{BrGo} 
if there exists a $\KK$-linear map $\partial \colon C \to \Der_\KK(R):z \mapsto \partial_z$ such that
\begin{enumerate}
\item[(i)] $C$ is stable under $\partial_z$ for all $z \in C$; and
\item[(ii)] the induced bracket $\{.,.\}$ on $C$ given by $\{z_1, z_2 \}:= \partial_{z_1}(z_2)$ turns $C$ into a Poisson algebra.
\end{enumerate}
For the definition and basic properties of the \emph{symplectic core} of a maximal ideal of an affine Poisson $\KK$-algebra $C$ as above, see \cite[\S 3.5]{BrGo}. Note
that symplectic cores are a generalisation of symplectic leaves (defined in the case when $\KK = \Cset$), 
and reduce to the latter when the Poisson bracket is \emph{complex algebraic},
\cite[Proposition 3.6]{BrGo}. 

In \cite[Theorem 4.2]{BrGo} it was proved that, if the base field is $\Cset$ and the Poisson order $R$ is $\Cset$-affine, then for all 
$\mm, \n \in \maxspec C$ in the same symplectic core of $\maxspec C$ with respect to the Poisson structure (ii), we have
\[
R/\mm R \cong R/\n R.
\]
This implies that

(*) if $(R, Z(R))$ is a $\Cset$-affine prime Poisson order, then $\AA(R)$ is a union of symplectic cores of $\maxspec Z(R)$
with respect to the Poisson structure (ii).
\\
In \cite[Theorem 3.10]{NTY} it was proved that

(**) if $(R, C)$ is a Poisson order over a base field of an arbitrary characteristic and 
$\tr \colon R \to C$ is a trace map that commutes with all derivations $\partial_z$ for $z \in C$, 
then $D_\ell(R/C, \tr)$ and $MD_\ell(R/C, \tr)$ are Poisson ideals of the Poisson algebra $C$ for all positive integers $\ell$.

Our Main Theorem thus provides a bridge between (*) and (**) because the zero locus of a Poisson ideal of a Poisson algebra $A$ is a union
of symplectic cores of $\maxspec A$.
\bigskip

\noindent
(b) \textbf{Normality of the center:} Given the importance of $Z(R)$ being normal in applying the Main Theorem, it is useful to have conditions on $R$ which ensure that
this property holds. We give two such conditions here, each of them sufficient to guarantee normality of the center, each of them frequently satisfied by PI algebras
occurring in representation theory.

(i) \textbf{Homological homogeneity:} Recall that a ring $R$ which is a finite module over its noetherian center is \emph{homologically homogeneous (hom. hom.)} if it
has finite global dimension $n$, and if, for every pair of simple $R$-modules $V$, $W$ with the same central annihilator, that is $\mathrm{Ann}_{Z(R)}(V)=
\mathrm{Ann}_{Z(R)}(W)$, $V$ and $W$ have the same projective dimension, $\mathrm{pr.dim}_R(V) = \mathrm{pr.dim}_R(W)$, \cite[$\S$1]{BrH}, see also
\cite[Definition 5.1]{BrMac}. Every hom. hom. ring is a finite direct sum of prime hom. hom. rings, \cite[Theorem 5.3]{BrH}. When $R$ is prime hom. hom. with $Z(R)$
affine over a field $\KK$, $\mathrm{pr.dim}_R(V)$ is in fact constant across \emph{all} simple $R$-modules, equalling $n$, \cite[Theorems 5.3, 4.8]{BrMac}.
\bpr{homnorm}\cite[Theorem 6.1]{BrH}
If $R$ is a prime hom. hom. ring, then $Z(R)$ is a normal noetherian domain.
\epr
\begin{proof} (Sketch) Let $P$ be a prime ideal of $R$ of height 1, and set $\mathfrak{p} := P \cap Z(R)$. Then $R_P := R \otimes_{Z(R)} Z(R)_{\mathfrak{p}}$ is hom.
hom. by \cite[Theorem 3.5]{BrH}. Since the global dimension of a hom. hom. ring equals its Krull dimension, by \cite[Theorem 2.5]{BrH}, $R_P$ is hereditary. Thus
$R_P$ is a hereditary noetherian prime ring satisfying a polynomial identity, so that $Z(R_P)$ is a Dedekind domain by \cite[Theorem 13.9.16]{McR}.

Now $R$ is a Cohen--Macaulay $Z(R)$-module by \cite[Theorem 2.5]{BrH}, so it follows \cite[Lemma 5.1(ii)]{BrH} that
\[
R = \cap \{R_P: \mathrm{height}(P) = 1 \}.
\]
Hence
\[
Z(R) = \cap \{Z(R_P): \mathrm{height}(P) = 1 \}.
\]
Thus $Z(R)$ is an intersection of normal subrings of its quotient field, and is therefore normal.
\end{proof}

It follows from results of Stafford-Zhang and Yi, \cite[Theorem 3.10]{StZ}, \cite[Proposition 3.2]{Y}, that an affine $\KK$-algebra $R$ which is a finite module over its
center is hom. hom. if and only if it is Auslander-regular and GK-Cohen--Macaulay. Bearing this in mind, one can assemble a large collection of ``standard'' examples of
hom. hom. algebras, including the following:
\begin{itemize}
\item[(1)]\cite[\S 6.2]{BZ} prime noetherian affine Hopf algebras of finite global dimension which satisfy a polynomial identity;
\item[(2)] (special case of (1)) enveloping algebras of finite dimensional restricted Lie algebras over a field of positive characteristic;
\item[(3)] (special case of (1)) noetherian prime regular group algebras satisfying a polynomial identity; that is, group algebras $\KK G$ of finitely generated abelian-by-
finite groups $G$ with no non-trivial finite normal subgroups and with no elements of order $\charr \KK$;
\item[(4)] (special case of (1)) quantized enveloping algebras and quantized function algebras with parameter $q$ a root of unity;
\item[(5)] (definition in \cite{EG}, sketch proof at \cite[$\S$4.4]{Br}) symplectic reflection algebras $H_{t,\mathbf{c}}$ with parameter $t = 0$.
\end{itemize}

(ii) \textbf{Maximal orders:} Maximal orders are defined at \cite[$\S$5.1.1]{McR}, with the standard reference being \cite{MR}. This class of rings constitutes a noncommutative analog of commutative normal domains, so it is not surprising that one easily shows \cite[Proposition 5.1.10(b)(i)]{McR} that
\[
\textit{the center of a prime Noetherian maximal order is normal.}
\]
Moreover, suppose that $R$ is a prime Noetherian ring satisfying a polynomial identity. If $Z(R)$ is normal, then $R$ equals its own trace ring $T(R)$,
\cite[13.8.2 and Proposition 13.8.11]{McR}, so that, if $R$ is in addition affine over a field, then $R$ is a finitely generated module over its center, by \cite[Proposition 13.9.11(ii)]{McR}.  Combining these observations with the classical results of Chamarie on
maximal orders cited below, we at once deduce that the following are prime Noetherian, affine $\KK$-algebras which are finite modules over their normal centers:
\begin{itemize}
\item[(1)]\cite[Proposition V.2.5]{MR} iterated skew polynomial $\KK$-algebras $R$, which satisfy a PI; that is, R satisfies a PI, and
$R = \KK \langle X_1, \dots X_n \rangle$, where, for all $i = 1, \ldots , n$,
\[
\KK\langle X_1, \ldots ,X_i\rangle = \KK\langle X_1, \ldots ,X_{i-1}\rangle[X_i; \sigma_i, \delta_i ],
\]
for $\KK$-algebra automorphisms $\sigma_i$ and $\sigma_i$-derivations $\delta_i$ of $\KK\langle X_1,  \ldots ,X_{i-1}\rangle$;
\item[(2)]\cite[Th$\acute{\mathrm{e}}$or$\grave{\mathrm{e}}$me X.2.1]{MR}, \cite[Theorem 5.1.6]{McR} $\KK$-algebras $R$ satisfying a polynomial identity,
and having an ascending $\mathbb{N}$-filtration $R = \cup_{i \geq 0}R_i$ with $R_0 = \KK$, whose associated graded $\KK$-algebra is a Noetherian prime maximal
order.
\end{itemize}

\noindent Note in passing that (ii)(2) yields an alternative way to incorporate the class (i)(2).
\sectionnew{An analog of the main theorem for Cayley-Hamilton algebras}
\label{CH}
\subsection{Statement of the theorem}
\label{statCH}
Let $R$ be a Cayley--Hamilton algebra of degree $n$ over an algebraically closed
field $\KK$ with trace $\tr : R \to C$, as in \S\ref{2.1extra}(4). Recall that $R$ is assumed to be
an affine algebra and the base field $\KK$
is assumed to have characteristic $\charr \KK \notin [1,n]$. In particular, $\tr(1) = n \neq 0$. Since the trace is
$C$-linear,
\[
C = \Im \tr.
\]
It follows that $C$ is also affine and that $R$ is a finitely generated $C$-module, see \cite[Theorem 4.5(a)]{DP}.

For $\mm \in \maxspec C$, denote by $\Irr_\mm(R)$ the equivalence classes of finite dimensional irreducible representations of $R$,
annihilated by $\mm$. The Azumaya locus of $R$ is defined as
\begin{multline*}
\AA(R) = \{ \mm \in \maxspec C \mid \mbox{$R$ has an irreducible representation} \\\mbox{of dimension $n$ annihilated by $\mm$} \}.
\end{multline*}
\bth{tCH} Let $R$ be a Cayley--Hamilton algebra of degree $n$ over an algebraically closed field $\KK$ of characteristic 0.
\begin{enumerate}
\item[(a)] The zero sets of the $n^2$-discriminant ideal and the modified $n^2$-discriminant ideal of $R$ over $C$ coincide
and are equal to the complement of the Azumaya locus $\AA(R) \subset \maxspec C$:
\[
\VV(D_{n^2}(R/C,\tr)) = \VV(MD_{n^2}(R/C,\tr)) = \maxspec C\; \backslash \;\AA(R).
\]
\item[(b)] For all positive integers $\ell$,
the zero sets of the $\ell$-discriminant ideal and the modified $\ell$-discriminant ideal of $R$ over $C$ coincide and are given in
terms of the irreducible representations of $R$ by
\begin{multline*}
 \VV(MD_{\ell}(R/C,\tr) )  = \VV(D_{\ell}(R/C,\tr))
\\
=  \Big\{ \mm \in \maxspec C \mid
\sum_{V \in \Irr_\mm(R)} (\dim_\KK V)^2  < \ell \Big\}.
\end{multline*}
\end{enumerate}
\eth
\subsection{Proofs}
\label{proofs}
As in Sect. \ref{main-sect}, for $\mm \in \maxspec C$ denote the canonical homomorphism
\[
\phi_\mm : R \to R/ \mm R
\]
and the trace map
\[
\tr_\mm : R /  \mm  R \to \KK \quad \mbox{given by} \quad
\tr_\mm (\phi_\mm(a)) := \phi_\mm(\tr(a)), \; \; a \in R.
\]

Extending the terminology of \S\ref{1.2}, we call a trace $\tr : R \to C$ representation theoretic if there exists a 
function $s_\mm : \Irr_\mm(R) \to \KK$ (where $\mm \in \maxspec C$) satisfying \eqref{rep-theor-tr2}.  

For every Cayley--Hamilton algebra $R$ of degree $n$ over an algebraically closed field $\KK$ of characteristic 0,
Procesi \cite[Theorem 2.6]{P}, \cite[Theorem 4.3]{DP} constructed a commutative algebra $\FF(R)$ with a $GL_n(\KK)$ action and
an embedding $\iota : R \hra M_n(\FF(R))$ of Cayley-Hamilton algebras such that
\[
\Im \iota = \FF(R)^{GL_n(\KK)} \quad \mbox{and} \quad \iota(C) = \FF(R)^{GL_n(\KK)}.
\]
We refer the reader to \cite[Sect. 1]{L} for a detailed exposition of this theorem and the
related background in invariant theory.

Recall that a representation $\phi : R \to \End_\KK(W)$ of $R$ is called {\em{trace preserving}} if
\[
\phi \circ \tr = \tr_W \circ \phi 
\]
where we identify $\KK \cong \KK \id_W$. 
For each $\mm \in \maxspec C$, Procesi's theorem produces a trace-preserving semisimple representation $W_\mm$ of $R/ \mm R$
such that each $V \in \Irr_\mm(R)$ is a subrepresentation of $W_\mm$, see \cite[Theorem 4.5(b)-(d) and Proposition 4.3]{DP}. Denoting 
by $k_\mm(V)$ the multuiplicity of $V \in \Irr_\mm(R)$ in $W_\mm$, gives the following:

\bth{P} [Procesi] Let $R$ be a Cayley--Hamilton algebra of degree $n$ over an algebraically closed field $\KK$ of characteristic 0.
For every $\mm \in \maxspec C$ there exists a function
\[
k_\mm : \Irr_\mm(R) \to \Zset_+
\]
such that
\[
n = \sum_{V \in \Irr_\mm(R)} k_\mm(V) \dim_\KK V
\]
and
\[
\tr_\mm = \sum_{V \in \Irr_\mm(R)} k_\mm(V) \tr_V
\]
where $\tr_V : \End_\KK(V) \to \KK$ is the usual trace.
\eth
Procesi's theorem proves that for every Cayley-Hamilton algebra
$R$ over a field of characteristic 0, its trace map $\tr R \to C$ is representation theoretic.
\begin{proof}[Proof of \thref{tCH}] (b) Let $\mm \in \maxspec C$. By an easy extension of \leref{2nd}, 
the trace $\tr_\mm : R/ \mm R \to \KK$ descends to a trace map
\[
\ol{\tr}_\mm : (R/ \mm R) / J( R/ \mm R) \to \KK.
\]
Furthermore, it follows from \thref{P} that the trace form
\[
\lcor -,  - \rcor_{\ol{\tr}_\mm}
\]
on $(R/ \mm R) / J( R/ \mm R)$ is non-degenerate. As in the proof of the Main Theorem,
\[
\dim_\KK \big( (R/ \mm R) / J( R/ \mm R) \big) = \sum_{V \in \Irr_\mm(R)} (\dim_\KK V)^2.
\]
If
\[
\ell \leq \sum_{V \in \Irr_\mm(R)} (\dim_\KK V)^2,
\]
then there exist $(y_1, \ldots, y_{\ell})$ and $(y'_1, \ldots, y'_{\ell}) \in R^{\ell}$ such that
\[
\phi_\mm \big( \det \big( [ \tr (y_i y'_j)]_{i,j=1}^{\ell} \big) \big) =
\det( [ \lcor \phi_\mm(y_i) \phi_\mm(y'_j) \rcor_{\ol{\tr}_\mm}]_{i,j=1}^{\ell}) \neq 0.
\]
If
\[
\ell > \sum_{V \in \Irr_\mm(R)} (\dim_\KK V)^2,
\]
then for all $(y_1, \ldots, y_{\ell})$ and $(y'_1, \ldots, y'_{\ell}) \in R^{\ell}$
\[
\phi_\mm \big( \det \big( [ \tr (y_i y'_j)]_{i,j=1}^{\ell} \big) \big) =
\det( [ \lcor \phi_\mm(y_i) \phi_\mm(y'_j) \rcor_{\ol{\tr}_\mm}]_{i,j=1}^{\ell}) = 0.
\]
The last two facts imply the statement of part (b) of the theorem.

(a) By \thref{P},
\[
\sum_{V \in \Irr_\mm(R)} (\dim_\KK V)^2 = n^2 \quad \mbox{if} \quad \mm \in \AA(R)
\]
and
\[
\sum_{V \in \Irr_\mm(R)} (\dim_\KK V)^2 < n^2 \quad \mbox{if} \quad \mm \in \maxspec C \backslash \AA(R).
\]
Combining these two facts with part (b) of the theorem for $\ell = n^2$, gives
\begin{multline*}
\VV(D_{n^2}(R/C,\tr)) = \VV(MD_{n^2}(R/C,\tr)) =  \\
\Big\{ \mm \in \maxspec C \mid
\sum_{V \in \Irr_\mm(R)} (\dim_\KK V)^2  < n^2 \Big\}
= \maxspec C\; \backslash \;\AA(R).
\end{multline*}
\end{proof}
\sectionnew{Singular loci and discriminant ideals}
\label{consequences}
\subsection{The singular locus of the center of a PI algebra}
\label{4.1}
Let $\mathcal{S}(C)$ denote the \emph{singular locus} of the commutative affine $\KK$-algebra $C$; that is,
\[
\mathcal{S}(C) = \{\mm \in \maxspec C : \mathrm{gl.dim} \, C_{\mm} < \infty \},
\]
where $\mathrm{gl.dim} \, U$ denotes the \emph{global dimension} of the ring $U$.
In this section we will
abbreviate
\[
D_{n^2}(R/Z(R)) := D_{n^2}(R/Z(R), \tr),
\]
suppressing the dependance of the discriminant on the trace.

\bco{gldim}Assume that $R$ and $\mathrm{tr}$ satisfy the hypotheses of the Main Theorem.
\begin{enumerate}
\item[(a)] Suppose that $R$ has finite global dimension. Then
\begin{equation}\label{include} \mathcal{S}(Z(R)) \subseteq \VV(D_{n^2}(R/Z(R))). \end{equation}
\item[(b)] Suppose that $R$ is a Cohen--Macaulay $Z(R)$-module, and that
\begin{equation}
\label{tall}
\mathrm{height}( D_{n^2}(R/Z(R))) \geq 2.
\end{equation}
Then
\begin{equation}
\label{reverse}
\VV(D_{n^2}(R/Z(R)))  \subseteq \mathcal{S}(Z(R)).
\end{equation}
\item[(c)] Suppose that $R$ has finite global dimension, that $R$ is a Cohen--Macaulay $Z(R)$-module, and that the inequality \eqref{tall} holds. Then
\begin{equation}
\label{equal}
\mathcal{S}(Z(R)) = \VV(D_{n^2}(R/Z(R))). \end{equation}
\end{enumerate}
\eco
\begin{proof}(a) Let $\mm \in \AA(R)$. Then $R_{\mm}$ is a free $Z(R)_{\mm}$-module. Now $\mathrm{gl.dim} \, R_{\mm} \leq \mathrm{gl.dim} \, R < \infty$. Therefore
$\mathrm{gl.dim} \, Z(R)_{\mm} $ is finite, since the global dimension of a commutative local noetherian ring is determined
by the projective dimension of its simple module,
and since a finite $R_{\mm}$-projective resolution of $R_{\mm}/\mm R_{\mm}$ affords a finite free $Z(R)_{\mm}$-resolution of a finite direct sum of copies of
$Z(R)_{\mm}/\mm Z(R)_{\mm}$. Therefore
\[
\mathcal{S}(Z(R)) \qquad \subseteq \qquad \maxspec R \;\setminus \;\mathcal{A} (R),
\]
and the result follows from the Main Theorem.

(b) The inequality \eqref{tall} coupled with the Main Theorem ensures that $R$ is Azumaya in codimension one. Thus all the hypotheses of \cite[Theorem 3.13]{BrMac}
are satisfied, so we can conclude from that result that every smooth point of $\maxspec Z(R)$ is Azumaya. Therefore the Main Theorem ensures that \eqref{reverse}
holds.

(c) This is immediate from (a) and (b).
\end{proof}

\subsection{Equivalent characterizations of the equality in \coref{gldim}(c)}
\label{4.2}
With a little more care, we can see that the hypothesis in \coref{gldim}(c), that the discriminant locus is ``small", is necessary:

\bco{sing} Assume that $R$ and $\mathrm{tr}$ satisfy the hypotheses of the Main Theorem. Suppose that $R$ has finite global dimension, that $R$ is a
Cohen--Macaulay $Z(R)$-module, and that $Z(R)$ is normal. Then the following are equivalent.
\begin{enumerate}
\item[(i)] $R$ is Azumaya in codimension one.
\item[(ii)] $\mathrm{height}( D_{n^2}(R/Z(R))) \geq 2.$
\item[(iii)] $\VV(D_{n^2}(R/Z(R))) = \mathcal{S}(Z(R)).$
\item[(iv)] $\VV(D_{n^2}(R/Z(R))) = \mathcal{S}(Z(R)) = \maxspec Z(R) \;\setminus \;\mathcal{A}(R).$
\item[(v)] $\VV(D_{n^2}(R/Z(R))) \subseteq \mathcal{S}(Z(R)).$
\end{enumerate}
\eco
\begin{proof} (i)$\Leftrightarrow$(ii): This follows from the Main Theorem.

(ii)$\Rightarrow$(iii): \coref{gldim}(c).

(iii)$\Rightarrow$(ii): This is immediate from the normality of $Z(R)$ and the fact that a normal noetherian domain of Krull dimension one is hereditary.

(iii)$\Leftrightarrow$(iv): One direction is trivial, and the other is given by the Main Theorem.

(iv)$\Rightarrow$(v): Trivial.

(iii)$\Rightarrow$(v): Clear.

(v)$\Rightarrow$(iii): \coref{gldim}(a).
\end{proof}
\sectionnew{Azumaya loci of PI quantized Weyl algebras}
\label{qWeyl}
\subsection{Setting}
\label{5.1}
The $n$-th quantized Weyl algebra for the parameters $E:=(\ep_1, \ldots, \ep_n) \in (\KK^\times)^n$
is the unital associative $\KK$-algebra $A^E_n$ with generators
\[
x_1, y_1, x_2, y_2, \ldots, x_n, y_n
\]
and relations
\begin{align*}
y_i y_k &= y_k y_i, \quad & \forall i, k,
\\
x_i x_k &= \ep_i x_k x_i, \quad &i< k,
\\
x_i y_k &= y_k x_i, \quad &i<k,
\\
x_i y_k &= \ep_k y_k x_i, \quad &i>k,
\\
x_i y_i - \ep_i y_i x_i &= 1 + \sum_{j=1}^{i-1} (\ep_j -1) y_j x_j, \quad &\forall i.
\end{align*}

Fix a partition $\la := (n_1 \geq \ldots \geq n_k > 0)$ of $N$ (so, $N= n_1 + \cdots + n_k$). Let $E := (\ep_1, \ldots, \ep_N) \in (\KK^\times)^N$.
Consider the $k$-fold tensor product
\begin{equation}
\label{A}
A^E_\la := A_{n_1}^{E_1} \otimes \cdots \otimes A_{n_k}^{E_k}
\end{equation}
where $E_l \in (\KK^\times)^l$ are the component subvectors of $E$ given by $E = (E_1, \ldots, E_k)$.

The algebra $A_\la^N$ is $\Zset^N$-graded
by setting
\[
\deg x_i = - \deg y_i := e_i
\]
where $\{e_1, \ldots, e_N \}$ is the standard basis of $\Zset^N$. Given a multiplicative antisymmetric bicharacter
\[
\chi \colon \Zset^N \times \Zset^N \to \KK^\times,
\]
denote by $A_{\la, \chi}^E$ the corresponding cocycle twist of $A_\la^E$. It is the unital associative $\KK$-algebra, identified with
$A_\la^E$ as a vector space, and with product
\[
a b = \chi(f,g) a * b \quad \mbox{for} \quad f, g \in \Zset^N, a \in (A_\la^E)_f, b \in (A_\la^E)_g
\]
where $a * b$ is the product of $A_\la^E$. The antisymmetric condition on $\chi$ means that
\[
\chi(f,g) \chi(g,f) =1, \; \; \chi(f,f) =1 \quad \mbox{for} \quad f, g \in \Zset^N.
\]

Given $i, l \in [1,N]$, set $i \prec l$ if
\[
n_1 + \cdots + n_{m-1} < i < l \leq n_1 + \cdots + n_m
\]
for some $m \in [1,k]$. Here and below the first sum is set to be 0 for $m=1$.
In the special case when the bicharacter $\chi$ satisfies
\[
\chi(e_i, e_l) = 1 \quad \mbox{for} \quad i \not\preceq l
\]
the algebra $A^E_{\la, \chi}$ is isomorphic to a product of mutiparameter quantized Weyl algebras,
and all such products arise in this way.
\subsection{Centers of PI quantized Weyl algebras}
\label{5.2}
For a root of unity $c \in \KK^\times$ denote its order
\[
\ord (c) := \min \{ m \in \Zset_+ : c^m = 1 \}.
\]

\bth{cent} Assume that $E$ and $\chi$ are taking values in roots of unity, and $\ep_i \neq 1$, $\forall i$.
The following are equivalent for the quantized Weyl algebra $A^E_{\la, \chi}$:
\begin{enumerate}
\item The algebra $A^E_{\la, \chi}$ is free over its center;
\item The center $\ZZ(A^E_{\la,\chi})$ is a polynomial algebra;
\item $\ord (\ep_i)$ and $\ord (\chi(e_i, e_j))$ divide $\ord (\ep_l)$ for all $i \preceq l$ and $j \in [1,N]$.
\end{enumerate}
If these conditions are satisfied, then
\[
\ZZ(A^E_{\la,\chi}) = \KK[x_i^{\ord (\ep_i)}, y_i^{\ord (\ep_i)}, 1 \leq i \leq N]
\]
and
\begin{equation}
\label{basis}
\{ x_1^{m_1} y_1^{m'_1} \ldots x_N^{m_N} y_N^{m'_N} : 1 \leq m_i, m'_i \leq \ord(\ep_i) \; \; \mbox{for} \; \; i \in [1,N]\}
\end{equation}
is a $\ZZ(A^E_{\la,\chi})$-basis of $A^E_{\la,\chi}$.
\eth
The theorem was proved in the case $k=1$ in \cite[Theorem A]{LY}. The proof of \cite{LY}
directly extends to the general case.

Let $i \in [1,N]$. Then
\begin{equation}
\label{i--m}
n_1 + \cdots + n_{m-1} < i \leq n_1 + \cdots + n_m \quad \mbox{for some} \quad m \in [1,k].
\end{equation}
Set
\[
z_i := [x_i,y_i] = 1 + (\ep_{n'} -1) y_{n'} x_{n'} + \cdots + (\ep_i -1) y_i x_i \in A^E_{\la, \chi}
\]
where
\begin{equation}
\label{nprime}
n' := n_1 + \cdots + n_{m-1} + 1.
\end{equation}
The elements $z_i$ are normal elements of $A^E_{\la, \chi}$:
\begin{equation}
\label{z-normal}
z_i x_l = \ep_l^{-\de_{l \preceq i}} x_l z_i , \quad
z_i y_l = \ep_l^{\de_{l \preceq i}} y_l z_i \quad \mbox{for} \quad i, l \in [1,N].
\end{equation}
\subsection{Azumaya loci and discriminant ideals of PI quantized Weyl algebras}
\label{5.3}
Assume that the three equivalent conditions in \thref{cent} are satisfied and set
\[
d_i := \ord (\ep_i) \quad \mbox{for} \quad i \in [1,N].
\]
Then
\[
\ZZ(A^E_{\la,\chi}) = \KK[X_i, Y_i, 1 \leq i \leq N]
\quad \mbox{where} \quad
X_i := x_i^{d_i}, Y_i := y_i^{d_i}.
\]
It follows from \eqref{z-normal} that the elements $Z_i := z_i^{d_i}$, $i \in [1,N]$ are also central. Let $i$ be as in \eqref{i--m}.
The element $Z_i$ is a polynomial in $\{X_{n'}, Y_{n'}, \ldots, X_i, Y_i\}$, recall \eqref{nprime}. It is
recursively given by
\begin{equation}
\label{Zgen}
Z_i = - (1 - \ep_i)^{d_i} Y_i X_i +
\begin{cases}
Z_{i-1}^{d_i/d_{i-1}}, & \mbox{if} \; \; i > n'  \\
1,          & \mbox{if} \; \; i = n'.
\end{cases}
\end{equation}
This follows from \cite[Eqs. (1.7) and (3.7)]{LY}.

\bth{qW-Azum-Discr} Assume that $E$ and $\chi$ are taking values in roots of unity, and $\ep_i \neq 1$, $\forall i$.
Assume also that the equivalent conditions in \thref{cent} are satisfied. Set $r:=d_1 \ldots d_N$.
\begin{enumerate}
\item[(a)] The discriminant of $A^E_{\la, \chi}$ over its center {\em{(}}with respect to the regular trace $\tr_{\reg}${\em{)}} is given by
\begin{align*}
d(A^E_{\la, \chi}/ \ZZ(A^E_{\la, \chi}), \tr_{\reg}) &=_{\KK^\times} r Z_1^{r^2 (d_1-1)/d_1} \ldots Z_N^{r^2 (d_N-1)/d_N}
\\
&=_{\KK^\times}  r z_1^{r^2 (d_1-1)} \ldots z_N^{r^2 (d_N-1)}
\end{align*}
provided that $\charr \KK \nmid r$.
\item[(b)] Assume that the base field $\KK$ has characteristic 0 or $\charr \KK > r^2$.
The algebra $A^E_{\la,\chi}$ has PI degree $r$ and its Azumaya locus is given by
\[
\{ \mm \in \maxspec \ZZ(A^E_{\la,\chi}) : Z_1(\mm) \neq 0, \ldots, Z_N(\mm) \neq 0 \}.
\]
\end{enumerate}
\eth
The special case of \thref{qW-Azum-Discr}(b) when $\la=(1,\ldots, 1)$, $\ep_1 = \ldots = \ep_N = \ep$ for a root of unity $\ep \neq 1$ and
$\chi(e_i,e_l) = \ep^{m_{il}}$ for some $m_{il} \in \Zset$ was previously obtained by Ganev \cite{Gan} and Cooney \cite{Coo}.
The special case of \thref{qW-Azum-Discr}(a) when $k=1$ was obtained in \cite{LY}.
\subsection{Proof of \thref{qW-Azum-Discr}}
\label{5.4}
(a) Using \cite[Theorem 8.2]{GY}, one can prove that $A^E_{\la,\chi}$ is isomorphic to a quantum cluster algebra when
all parameters $\ep_i \in \KK^\times$ are not roots of unity. This structure can then be specialized to the roots of unity case.
We will only need two of its clusters which we construct directly. Denote the skewpolynomial algebras
\begin{align*}
&S_x :=
\frac{\KK\lcor x_i , z_i, \; \; 1 \leq  i \leq N \rcor}
{(x_i x_l - \ep_i^{\de_{i \preceq l}} \chi(e_i, e_l)^2 x_l x_i,  \; \; z_i z_l = z_l z_i, \;  \; z_i x_l = \ep_l^{-\de_{l \preceq j}} x_l z_i)} \; \; \mbox{and}
\\
&S_y :=
\frac{\KK\lcor y_i, z_i, \; \; 1 \leq  i \leq N \rcor}
{(y_i y_l = \chi(e_l,e_i)^2 y_l y_i, \; \; z_i z_l = z_l z_i, \; \; z_i y_l = \ep_l^{\de_{l \preceq i}} y_l z_i)} \cdot
\end{align*}
We have the isomorphisms
\begin{align}
\label{iso1}
&A^E_{\la, \chi}[x_i^{-d_i}, 1 \leq i \leq N ] \cong S_x [x_i^{-d_i}, 1 \leq i \leq N] \; \;
\mbox{and} \quad
\\
&A^E_{\la, \chi}[y_i^{-d_i}, 1 \leq i \leq N] \cong S_y [y_i^{-d_i}, 1 \leq i \leq N].
\label{iso2}
\end{align}
Indeed, it follows from \eqref{z-normal} and the definition of the elements $z_i$ that the algebra $A^E_{\la, \chi}$ embeds in the right hand sides of
\eqref{iso1}--\eqref{iso2}. Furthermore, the target algebras are generated by the images of $A^E_{\la, \chi}$, $\{x_i^{-d_i} : 1 \leq i \leq N\}$
and $\{y_i^{-d_i} : 1 \leq i \leq N\}$, respectively. This proves \eqref{iso1}--\eqref{iso2}. Clearly,
\[
\ZZ(S_x) = \KK[x_i^{d_i}, z_i^{d_i}, 1 \leq i \leq N] \quad \mbox{and} \quad
\ZZ(S_y) = \KK[y_i^{d_i}, z_i^{d_i}, 1 \leq i \leq N].
\]
Moreover $S_x$ and $S_y$ are free modules over $\ZZ(S_x)$ and $\ZZ(S_y)$ with bases
\begin{equation}
\{ x_1^{m_1} z_1^{m'_1} \ldots x_N^{m_N} z_N^{m'_N} : 1 \leq m_i, m'_i \leq d_i \}, \quad
\{ y_1^{m_1} z_1^{m'_1} \ldots y_N^{m_N} z_N^{m'_N} : 1 \leq m_i, m'_i \leq d_i \},
\label{P-bases}
\end{equation}
respectively.

Denote the algebras in \eqref{iso1}--\eqref{iso2}  by $L_x$ and $L_y$. The extensions of the regular trace of $A^E_{\la, \chi}$ to
the localizations in $L_x$ and $L_y$ match the extensions of the regular traces of $S_x$ and $S_y$ to
$L_x$ and $L_y$. (In both case these extensions are simply the regular traces of $L_x$ and $L_y$ because the localizations
are central.) Therefore,
\begin{align}
\label{eqq1}
&d(A^E_{\la, \chi}/ \ZZ(A^E_{\la,\chi}), \tr_{\reg}) =_{L_x^\times} d (L_x / \ZZ(L_x), \tr_{\reg}) =_{L_x^\times} d(S_x/\ZZ(S_x), \tr_{\reg}),
\\
\label{eqq2}
&d(A^E_{\la,\chi}/ \ZZ(A^E_{\la,\chi}), \tr_{\reg}) =_{L_y^\times} d(L_y / \ZZ(L_y),  \tr_{\reg}) =_{L_y^\times} d(S_y/\ZZ(S_y),  \tr_{\reg}).
\end{align}
Since $S_x$ and $S_y$ are skewpolynomial algebras with bases \eqref{P-bases} over their centers,
their discriminats are given by \cite[Proposition 2.8]{CPWZ}
\begin{align}
\label{eqq3}
&d(S_x/\ZZ(S_x), \tr_{\reg}) =_{\KK^\times} r \prod_{i=1}^N x_i^{r^2 (d_i-1)} z_i^{r^2(d_i-1)}, \\
&d(S_y/\ZZ(S_y), \tr_{\reg}) =_{\KK^\times} r \prod_{i=1}^N y_i^{r^2 (d_1-1)} z_i^{r^2 (d_i-1)}.
\label{eqq4}
\end{align}
The groups of unites of $L_x$ and $L_y$ are the Laurent monomials in the variables $x_i$ and $y_i$, respectively.
It follows from \eqref{eqq1}--\eqref{eqq4} that
\begin{align*}
d(A^E_{\la, \chi}/ \ZZ(A^E_{\la,\chi}),  \tr_{\reg}) &=_{\KK^\times}  r
x_1^{m_1 d_1} \ldots x_N^{m_N d_N} z_1^{r^2 (d_1-1)} \ldots z_N^{r^2 (d_N-1)}
\\
&=_{\KK^\times} r y_1^{l_1 d_1} \ldots y_N^{l_N d_N} z_1^{r^2 (d_1-1)} \ldots z_N^{r^2 (d_N-1)}
\end{align*}
for some $m_i, l_i \in \Zset$. If $\charr \KK \mid r$, we are done. Otherwise,
this implies that
\[
x_1^{m_1 d_1} \ldots x_N^{m_N d_N} =_{\KK^\times} y_1^{l_1 d_1} \ldots y_N^{l_N d_N}
\]
because $A^E_{\la,\chi}$ is a domain. Since $\{ x_1^{m_1} \ldots x_N^{m_N} y_1^{l_1} \ldots x_N^{l_N} : m_i, l_i \in \Nset \}$ is a basis
of $A^E_{\la,\chi}$, $m_i = l_i =0$ for all $i$, which completes the proof of part (a).

(b) It follows from \eqref{basis} that $A^E_{\la, \chi}$ embeds in a matrix algebra of size $r \times r$ over
$\ZZ(A^E_{\la.\chi})$. Thus, the PI degree of $A^E_{\la, \chi}$ is at most $r$. Part (b) of the Main Theorem and the
non-vanishing of the discriminant in part (a) imply that the PI degree of $A^E_{\la,\chi}$ is not strictly smaller than $r$.

The statement for the Azumaya locus of $A^E_{\la, \chi}$ follows from the Main Theorem, part (a) of this theorem, and the fact that
\[
D_{r^2}(A^E_{\la,\chi}/\ZZ(A^E_{\la,\chi})) = \lcor d(A^E_{\la,\chi}/\ZZ(A^E_{\la,\chi}), \tr_{\reg}) \rcor
\]
coming from the freeness of $A^E_{\la,\chi}$ over its center, cf. \leref{poly}(b). This completes the proof of \thref{qW-Azum-Discr}.

\end{document}